\documentclass[onecolumn,amsmath,amssymb]{revtex4}

\usepackage[dvipdfmx]{graphicx}
\usepackage{dcolumn}
\usepackage{bm,color}
\usepackage{amsmath}
\usepackage{makeidx}

\newtheorem{theorem}{Theorem}
\newtheorem{corollary}[theorem]{Corollary}
\newtheorem{lemma}[theorem]{Lemma}

\newtheorem{definition}{Definition}

\makeindex
\begin{document}

\preprint{APS/}

\title
{Theory of $B(X)$-module \\
- Algebraic module structure of generally-unbounded infinitesimal generators -}

\author{Yoritaka Iwata}
 \email{iwata$\_$phys@08.alumni.u-tokyo.ac.jp}

 \affiliation{Faculty of Chemistry ,Materials and Bioengineering, Kansai University, Osaka 564-8680, Japan}

\begin{abstract}
The concept of logarithmic representation of infinitesimal generators is introduced, and it is applied to clarify the algebraic structure of  bounded and unbounded infinitesimal generators.
In particular, by means of the logarithmic representation, the bounded components can be extracted from generally-unbounded infinitesimal generators.
In conclusion the concept of module over a Banach algebra is proposed as the generalization of Banach algebra.
As an application to mathematical physics, the rigorous formulation of rotation group, which consists of unbounded operators being written by differential operators, is provided using the module over a Banach algebra.
\end{abstract}

\maketitle

\tableofcontents

\newpage

\section{Introduction} 
Based on the logarithmic representation of infinitesimal generators, a module over a Banach algebra is introduced.
Let us call such an algebraic subject the $B(X)$-module, where $X$ and $B(X)$ stand for a Banach space and its operator algebra, respectively.
The $B(X)$-module does not correspond only to the extension of Banach algebra, but also to the general authorization of the Lie algebra consisting of differential operators.
This algebraic entity is an operator algebra being introduced based on the framework of logarithmic representation of operators.
There are two concepts, which are to be bridged in this paper: a set of infinitesimal generators generating groups or semigroups of operators, and the elements of Lie algebra.
The following statements are valid: 
\begin{itemize}
\item the sum of two closed operators are not necessarily a closed operator, so that  the sum of two infinitesimal generators are not necessarily an infinitesimal generator;
\item the sum of two elements in the Lie algebra are necessarily an element of the Lie algebra.
\end{itemize}
Here is a contradiction in some general situations, as seen in the relation between the Lie group and the Lie algebra in which the Lie algebra corresponds to a set of infinitesimal generators.
Besides, these two statements are true if the two infinitesimal generators are bounded operators. 
More substantially, the product cannot be justified without limiting ourselves to (sub)sets of bounded operators.
In this paper, by means of the logarithmic representation of infinitesimal generators of invertible evolution operators, a set of generally-unbounded infinitesimal generators is characterized as an algebraic module over a Banach algebra.
The logarithm of operators is a key to make a bridge for these contradicting statements.

The logarithm of an injective sectorial operator was introduced by Nollau~\cite{69nollau} in 1969.
After a long time, the logarithm of sectorial operators were studied again from 1990's \cite{94boyadzhiev,00okazawa-1, 00okazawa-2}, and its utility was established with respect to the definition of the logarithms of operators~\cite{03hasse,01martinez} (for a review of sectorial operators, see Hasse \cite{06hasse}).
While the sectorial operator has been a generic framework to define the logarithm of operators, the sectorial property is not generally satisfied by the evolution operators.
In this sense, it is necessary to introduce a reasonable framework for defining the logarithm of non-sectorial infinitesimal generators.

In this paper, the theory of $B(X)$-module is introduced.
The utility of the theory is confirmed in the application to
\begin{enumerate}
\item the solvability of abstract  Cauchy problem
\item the generalization of the Cole-Hopf transform
\item the foundation of the rotation group
\end{enumerate}
in the latter parts of Sections II, III and IV.

This paper is the completion of the recent studies shown in Refs.~\cite{17iwata-1,17iwata-3,17iwata-2,18iwata-1,18iwata-2,19iwata-book,19iwata,20iwata}. 
By organizing the preceding works into a logical order, the several statements are renewed. 
First, the generalized version of logarithmic representation (Col. 5) is possible using the concept of alternative infinitesimal generator.
Even without any additional assumptions, it enables us to generalize the logarithmic representation for infinitesimal generators of non-invertible evolution operators.  
Although this fact is taken for granted in the lately published papers of Refs.~\cite{17iwata-1,17iwata-3,17iwata-2,18iwata-1,18iwata-2,19iwata-book,19iwata,20iwata}, it is mentioned within a logical process for the first time. 
Second, although the relativistic formulation is introduced for changing the evolution direction as seen in the application of the Cole-Hopf transform, it should not be restricted to the application of the Cole-Hopf transform.
The relativistic formulation of abstract evolution equation is a kind of generalization of abstract evolution equations. 
More clearly it generalize the concept of abstract evolution equation to the abstract equation.
Consequently the theory of $B(X)$-module is written in the relativistic form.

\section{Logarithmic representation of operators} 
\subsection{Banach algebra}
Let $(A,\| \cdot \|)$ be a Banach space\index{Banach space} (for a textbook, see \cite{65yosida}).
A mapping
\[ A \times A \to A, \quad (x,y) \to x \cdot y \]
is called a multiplication on $A$, if it is bilinear and associative.
$\| \cdot \|$ is said to be a submultiplicative norm if $ \| x \cdot y \| \le \| x\| ~\| y \|$ for  each $x, y \in A$.
The Banach space $A$ together with a multiplication and submultiplicative norm is called a Banach algebra\index{Banach algebra}.

Let X be a Banach space.
Denote by $B(X)$ the set of all bounded linear operators $f:X \to X$.
Then $A = B(X)$ is an example of a Banach algebra with multiplication as composition and norm defined by
\[ \| f \| = \sup \{ \|f(x)\| : \|x \| \le 1 \}, \]
where $f \in B(X)$.
$B(X)$ is called the operator algebra of $X$.
Let $A$ be a Banach algebra and let $X$ be a Banach space.
$X$ is said to be
\begin{enumerate}
\item A left Banach $A$-module\index{left Banach $A$-module} if there exists a bilinear mapping $\cdot: A \times X \to X$, $(a, x) \to a \cdot x$, called left module action, such that $\|a \cdot x \| \le \|a\| \|x \|$ and
\[ (ab)\cdot x = a \cdot (b \cdot x), \quad a, b \in A ~{\rm and}~ x \in X. \]
\item A right Banach $A$-module\index{right Banach $A$-module} if there exists a bilinear mapping $\cdot: A \times X \to X$, $(x, a) \to x \cdot a$, called right module action, such that $\|x \cdot a \| \le \|a\| \|x \|$ and
\[ x \cdot (ab) = (x \cdot a) \cdot b, \quad a, b \in A ~{\rm and}~ x \in X. \]
\item A Banach $A$-module\index{Banach $A$-module} if it is a left and right Banach $A$-module and
\[ a\cdot (x \cdot b) = (a \cdot x) \cdot b, \quad a, b \in A ~{\rm and}~ x \in X. \]
\end{enumerate}
As an example of Banach algebra, $B(X)$ is taken in the following.
Then $X$ is a Banach $B(X)$-module under 
\[ \begin{array}{ll}
B(X) \times X \to X, \quad (f, x) \to f(x),   \vspace{2.5mm}  \\
X \times B(X) \to X, \quad (x, f) \to f(x).
\end{array} \]
Banach space, Banach algebra, and Banach $B(X)$-module ($B(X)$-module, for short) are the basic concepts in this paper.

\subsection{Two parameter group on Banach spaces} \label{tp-group}
All the discussion begins with the definition of groups on the Banach spaces that will be generalized to a well-defined semigroup in later sections.
Let $(X,\| \cdot \|)$ be a Banach space and $B(X)$ its Banach algebra\index{operator algebra}.
In particular $B(X)$ is an example of a Banach algebra.

A two-parameter group\index{two-parameter group} is defined on $X$.
Let $X$ be a Banach space and $T \in (0, \infty)$.
A two parameter group on X is an operator valued mapping $(t, s) \to U(t, s)$ from $[-T,T]$ into $B(X)$ with the semigroup properties\index{semigroup property}:
\begin{equation} \begin{array}{ll}  \label{sg1}
U (t, r)U (r, s) = U (t, s), \qquad r, s, t \in [-T, T], \vspace{2.5mm}  \\
U (s, s) = I , \qquad  s \in [-T, T],  
\end{array} \end{equation}
and the strong continuity; for each $s \in [-T, T ]$ and $x \in X$, the map $t \to U (t, s)x$ is continuous on $[s,T]$.
Both $U(t,s)$ and $U(s,t)$ are assumed to be well-defined to satisfy 
\begin{equation} \label{sg3}
U(s,t) ~ U(t,s) = U(s,s) = I, 
\end{equation}
where $U(s,t)$ corresponds to the inverse operator of $U(t,s)$. 
Since $U(t,s) ~ U(s,t) = U(t,t) = I$ is also true, the commutation between $U(t,s)$ and $U(s,t)$ follows. 
Operator $U(t,s)$, which is called the evolution operator\index{evolution operator} in the following, is a generalization of exponential function; indeed the properties shown in Eqs.~\eqref{sg1}-\eqref{sg3} are satisfied by taking $U(t,s)$ as $e^{t-s}$.
Evolution operator is an abstract concept of exponential function valid for both finite and infinite dimensional Banach spaces.
Due to the validity of Eq.~\eqref{sg3}, the invertible evolution family is to be associated with some linear evolution equations of hyperbolic type and those of dispersive type. 
In the same context, the obtained results can be directly applied to some semilinear evolution equations\index{semilinear evolution equations}~(for a text book, see \cite{98cazenave}).
For example, the solutions of linear and nonlinear wave equations are written by the evolution operator $U(t,s)$ defined above. 

Let $Y$ be a dense Banach subspace of the Banach space $X$, and the topology of $Y$ be stronger than that of $X$.
The space $Y \subset X$ is assumed to be $U(t,s)$-invariant; for any $t,s$ satisfying $-T \le t, s  \le T$, $U(t,s) Y = Y$.
Following the definition of $C_0$-(semi)group (cf. the assumption $H_2$ in Sec.~5.3 of Pazy~\cite{83pazy} or corresponding discussion in Kato~\cite{70kato,73kato}), $U(t,s)$ trivially satisfy the boundedness in the present setting; there exist real numbers $M$ and $\beta$ such that
\begin{equation} \label{qb} 
\| U(t,s)  \|_{B(X)} \le M e^{\beta t},
\quad
\| U(t,s)  \|_{B(Y)} \le M e^{\beta t}.
\end{equation}
that are practically reduced to
\[ \| U(t,s)  \|_{B(X)} \le M e^{\beta T}, \quad \| U(t,s)  \|_{B(Y)} \le M e^{\beta T},  \]
when the interval is restricted to be finite $[-T,T]$.
Since $C_0$ semigroup theory is essentially based on the Laplace transform of operators, the satisfaction of Eq.~(\ref{qb}) is discussed here; $M e^{\beta t}$ in Eq. (\ref{qb}) arises from the condition for the existence theorem for the Laplace transforms (for example, see \cite{11kreyszig}), and $M' = M e^{\beta T}$ is regarded as a finite real number in the present setting.

Next, for the well-defined $U(t,s)$, the counterpart of the logarithm in the abstract framework is introduced.
There are two concepts associated with the logarithm of operators; one is the infinitesimal generator\index{infinitesimal generator} and the other is $t$-differential of $U(t,s)$.
These two concepts are connected as follows.

\begin{definition} [Pre-infinitesimal generator]
For $-T \le t, s  \le T$, the weak limit\index{weak limit}
\[ \begin{array}{ll}
 \mathop{\rm wlim}\limits_{h \to 0} h^{-1} (U(t+h,s) - U(t,s)) ~u_s 
= \mathop{\rm wlim}\limits_{h \to 0}   h^{-1}(U(t+h,t) - I) ~ U(t,s) ~u_s,  
\end{array} \]
is assumed to exist for certain $u_s$, which is an element of a dense subspace $Y$ of $X$. 
A linear operator $A(t): Y  ~\to~  X$ is defined by
\begin{equation} \label{pe-group}
A(t) u_t := \mathop{\rm wlim}\limits_{h \to 0}  h^{-1} (U(t+h,t) - I) u_t
\end{equation}
for $u_t \in Y$ and $-T \le t, s  \le T$.
The operator $A(t)$ for a whole family $\{U(t,s)\}_{-T \le t,s \le T}$ is called the pre-infinitesimal generator\index{pre-infinitesimal generator}. 
 \end{definition}

Let $t$-differential of $U(t,s)$ in a weak sense~\cite{19iwata-book} be denoted by
\begin{equation} \label{de-group} \begin{array}{ll} 
\partial_t U(t,s)~u_s = A(t) U(t,s) ~u_s.
 \end{array}  \end{equation}
Equation~\eqref{de-group} is regarded as a differential equation
satisfied by $u(t) = U(t,s) u_s$ that implies a relation between $A(t)$ and the logarithm:
\[ \begin{array}{ll}
A(t)  = [\partial_t U(t,s)] ~ U(s,t).
\end{array} \]
The relation between $A(t)$ and the logarithm is discussed in the next section \ref{sec1b}.
Pre-infinitesimal generators are not necessarily infinitesimal generators without assuming a dense property of domain space $Y$ in $X$.
For example, in $t$-independent cases, an operator $A(t)$ defined by  Eq.~\eqref{pe-group} is not necessarily a densely-defined and closed linear operator, while $A(t)$ must be a densely-defined and closed linear operator with its
resolvent set included in $\{\lambda \in {\mathbb C}: {\rm Re} \lambda > \beta
\}$ for $A(t)$ to be the infinitesimal generator.
On the other hand, infinitesimal generators are necessarily pre-infinitesimal generators.
That is, only the exponentiability with a certain ideal domain is valid to the pre-infinitesimal generators.
The definition of pre-infinitesimal generator is useful in terms of providing the algebraic structure. 
Let a set of pre-infinitesimal generators be denoted by $G(X)$.
It is trivial that $B(X) \subset G(X)$.

\subsection{Logarithmic representation of pre-infinitesimal generator} \label{sec1b}
The logarithmic representation of infinitesimal generator\index{logarithmic representation of infinitesimal generator} is introduced in order to clarify the structure of infinitesimal generators~\cite{17iwata-1}.
The logarithm of $U(t,s)$ is defined by the Dunford-Riesz integral~\cite{43dunford}.
The boundedness of $U(t,s)$ on $X$ makes the problem rather easy.
Indeed the boundedness allows us to introduce the translation on the complex plane as a tool to realize the parallel displacement of the entire spectral set.
On the other hand, two difficulties inherent to the logarithm
\begin{itemize}
\item singularity of logarithm at the origin
\item multi-valued property of the logarithm
\end{itemize}
arise.
By introducing a constant $\kappa \in {\mathbb C}$, the singularity can be handled.
This simple treatment is definitely practical to well-define the logarithm of non-sectorial operators.
By introducing a principal branch\index{principal branch}  (denoted by ``Log") of the logarithm (denoted by ``$ \log$"), the multi-valued property is handled.
Indeed, for any complex number $z \in C$, a  branch of logarithm is defined by
\[ \begin{array}{ll}
 {\rm Log} z  = \log |z| + i \arg Z,  
\end{array} \]
where $Z$ is a complex number chosen to satisfy $|Z| = |z|$, $-\pi < \arg Z \le \pi$, and $\arg Z = \arg z + 2 n \pi$ for a certain integer $n$. 

\begin{lemma} [Logarithmic representation of operators]   \label{lem3}
Let $t$ and $s$ satisfy $0 \le t,s \le T$.
For a given $U(t,s)$ defined in Sec.~\ref{tp-group}, its logarithm is well defined; there exists a certain complex number $\kappa$ satisfying
\begin{equation}
\label{logex13} \begin{array}{ll}
{\rm Log} (U(t,s)+\kappa I) = \frac{1}{2 \pi i} \int_{\Gamma} {\rm Log} \lambda  
 ~ ( \lambda I - U(t,s) - \kappa I )^{-1}  d \lambda,
\end{array} \end{equation}
where an integral path $\Gamma$, which excludes the origin, is a circle
in the resolvent set of $U(t,s) +\kappa I$.
Here $\Gamma$ is independent of $t$ and $s$. 
${\rm Log} (U(t,s)+ \kappa I)$ is bounded on $X$.
\end{lemma}

 {\bf Proof.}  
The logarithm ${\rm Log}$ holds the singularity at the origin, so that it is necessary to show a possibility of taking a simple closed curve (integral path) excluding the origin in order to define the logarithm by means of the Dunford-Riesz integral. 
It is not generally possible to take such a path in case of $\kappa =0$.

First, $U(t,s)$ is assumed to be bounded for $0 \le t,s \le T$ (Eq.~\eqref{qb}), and the spectral set of $U(t,s)$ is a bounded set in ${\mathbb C}$.
Second, for $\kappa$ satisfying
\[ |\kappa| > M e^{\beta T}, \]
the spectral set of $U(t,s) + \kappa I$ is separated with the origin.
Consequently it is possible to take an integral path $\Gamma$ including the spectral set of $U(t,s)+\kappa I$ and excluding the origin.
Equation~\eqref{logex13} follows from the Dunford-Riesz integral.
Furthermore, by adjusting the amplitude of $\kappa$, an appropriate integral path always exists independent of $t$ and $s$.
${\rm Log} (U(t,s)+\kappa I)$ is bounded on $X$, since $\Gamma$ is included in the resolvent set of $(U(t,s)+\kappa I)$.
\quad  $\square$   \\

According to this lemma, by introducing nonzero $\kappa$, the logarithm
of $U(t,s)+\kappa I$ is well-defined without assuming the sectorial property to $U(t,s)$.
On the other hand Eq.~\eqref{logex13} is valid with $\kappa=0$ only for limited cases.

\begin{theorem} [Logarithmic representation of infinitesimal generators]
\label{thm1}
Let $t$ and $s$ satisfy $-T \le t,s \le T$, and $Y$ be a dense subspace of $X$.
For $U(t,s)$ defined in Sec.~\ref{tp-group}, let $A(t) \in G(X)$ and $\partial_t U(t,s)$ be determined by Eqs.~\eqref{pe-group} and \eqref{de-group} respectively.  
If $A(t)$ and $U(t,s)$ commute, pre-infinitesimal generators $\{ A(t) \}_{-T \le t \le T}$ are represented by means of the logarithm function; there exists a certain complex number $\kappa \ne 0$ such that
\begin{equation} \label{logex} \begin{array}{ll}
 A(t) ~ u_s =  (I+ \kappa U(s,t))~ \partial_{t} {\rm Log} ~ (U(t,s) + \kappa I) ~ u_s, 
\end{array} \end{equation}
where $u_s$ is an element in $Y$.
Note that $U(t,s)$ defined in Sec.~\ref{tp-group} is assumed to be invertible.
\end{theorem}

{\bf Proof. ~ }
For $U(t,s)$ defined in Sec.~\ref{tp-group}, operators $ {\rm Log} ~ (U(t,s) + \kappa I)$ and $ {\rm Log} ~ (U(t+h,s) + \kappa I)$ are well defined for a certain $\kappa$ (Lemma~\ref{lem3}).
The $t$-differential in a weak sense\index{t-differential in a weak sense} is formally written by
\begin{equation} \label{difference0} \begin{array} {ll} 
\mathop{\rm wlim}\limits_{h \to 0}  \frac{1}{h} \{ {\rm Log} ~(U(t+h,s)+\kappa I) - {\rm Log} ~(U(t,s)+ \kappa I) \}   \vspace{1.5mm} \\
  =\mathop{\rm wlim}\limits_{h \to 0}   \frac{1}{h} \frac{1}{2 \pi i}
 \int_{\Gamma} {\rm Log} \lambda   
 ~ \{ ( \lambda I - U(t+h,s) - \kappa I )^{-1}   
  -  ( \lambda I - U(t,s) - \kappa I )^{-1}  \} d \lambda  \vspace{1.5mm} \\  
 = \mathop{\rm wlim}\limits_{h \to 0}  \frac{1}{2 \pi i}
 \int_{\Gamma} {\rm Log} \lambda   
 ~ \{ (\lambda I - U(t+h,s)-\kappa I )^{-1} \frac{U(t+h,s)-U(t,s)}{h} (\lambda - U(t,s)- \kappa I )^{-1} \} d \lambda   
\end{array} \end{equation} 
where $\Gamma$, which is possible to be taken independent of $t$, $s$ and $h$ for a sufficiently large certain $\kappa$, denotes a circle in the resolvent set of both $U(t,s)+ \kappa I$ and $U(t+h,s)+\kappa I$.
A part of the integrand of Eq.~\eqref{difference0} is estimated as
\begin{equation} \label{intee} \begin{array} {ll}
\quad \| \{ (\lambda I - U(t+h,s)-\kappa I )^{-1} \frac{U(t+h,s)-U(t,s)}{h} (\lambda I  - U(t,s)- \kappa I )^{-1} \} v \|_X \vspace{1.5mm} \\
  \le  \| \{ (\lambda I - U(t+h,s)- \kappa I )^{-1} \|_{B(X)}  
 \|  \frac{U(t+h,s)-U(t,s)}{h} (\lambda I - U(t,s)- \kappa I )^{-1} \} v \|_X, 
\end{array} \end{equation}
for $v \in X$. 
There are two steps to prove the validity of Eq.~\eqref{difference0}. \\

{\bf [1st step] ~}
the former part of the right hand side of Eq.~\eqref{intee} satisfies
\[ \begin{array} {ll}
 \| \{ (\lambda I - U(t+h,s)- \kappa I )^{-1} \|_{B(X)}  < \infty,
\end{array} \]
since $\lambda$ is taken from the resolvent set of $U(t+h,s)- \kappa I$.
In the same way the operator $(\lambda - U(t,s) - \kappa )^{-1}$ is bounded on $X$ and $Y$.
Then the continuity of the mapping $t \to (\lambda - U(t,s)- \kappa )^{-1}$ as for the strong topology follows:
\[ \begin{array} {ll}
 \|  (\lambda I - U(t+h,s)- \kappa I )^{-1} - (\lambda I - U(t,s)- \kappa I)^{-1}  \|_{B(X)}   \vspace{1.5mm} \\
 \le  \| (\lambda I - U(t+h,s)- \kappa I )^{-1} \|_{B(X)}   
 \|( U(t+h,s)-U(t,s)) (\lambda I - U(t,s)-\kappa I )^{-1}  \|_{B(X)}. \\
\end{array} \] 

{\bf [2nd step] ~}
The latter part of the right hand side of Eq.~\eqref{intee} is estimated as
\begin{equation} \label{unibound}  \begin{array}{ll}
 \left\| \frac{U(t+h,s)-U(t,s)}{h} (\lambda I - U(t,s)-\kappa I )^{-1}  u  \right\|_X \vspace{1.5mm}\\
 = \left\| \frac{1}{h} \int_t^{t+h} A(\tau) U(\tau,s) (\lambda I  - U(t,s)- \kappa I )^{-1}  u ~ d\tau  \right\|_X \vspace{1.5mm}\\
 \le  \frac{1}{|h|}  \int_t^{t+h} \| A(\tau)  U(\tau,s)\|_{B(Y,X)}  
   \| (\lambda I - U(t,s)- \kappa I )^{-1}\|_{B(Y)} \|u \|_{Y} ~ d\tau 
\end{array} \end{equation}
for $u \in Y$.
Because $\| A(\tau) U(\tau,s) \|_{B(Y,X)}  < \infty $ is true by assumption, the right hand side of Eq.~\eqref{unibound} is finite.  
Equation~\eqref{unibound} shows the uniform boundedness with respect to $h$, then the uniform convergence ($h \to 0$) of Eq.~\eqref{difference0} follows.
Consequently the weak limit process $h \to 0$ for the integrand of Eq. (\ref{difference0}) is justified, as well as the commutation between the limit and the integral.

According to Eq.~\eqref{difference0}, interchange of the limit with the integral leads to
\[ \begin{array} {ll} 
\partial_t {\rm Log} (U(t,s) + \kappa I) ~ u 
=  \frac{1}{2 \pi i} \int_{\Gamma} d\lambda  
 ~ \left[ ( {\rm Log} \lambda )  (\lambda I - U(t,s)- \kappa I )^{-1} 
 ~\mathop{\rm wlim}\limits_{h \to 0} \left ( \frac{U(t+h,s)-U(t,s)}{h} \right)
~ (\lambda I - U(t,s)- \kappa I)^{-1}    \right] ~u \\
\end{array} \] 
for $u \in Y$.
Because it is also allowed to interchange $A(t)$ with $U(t,s)$,
\[  \begin{array}{ll}
\partial_t {\rm Log} (U(t,s) + \kappa I) ~ u \vspace{1.5mm}\\
 = \frac{1}{2 \pi i} \int_{\Gamma} ({\rm Log} \lambda)
 (\lambda I -U(t,s)-\kappa I )^{-1} 
  A(t) ~ U(t,s) ~ (\lambda I -U(t,s)- \kappa I)^{-1}  d \lambda ~ u \vspace{1.5mm}\\
 = \frac{1}{2 \pi i} \int_{\Gamma} ({\rm Log} \lambda) ~(\lambda I -U(t,s) - \kappa I)^{-2} ~ U(t,s) ~ d \lambda ~A(t) ~ u
\end{array} \]
for $u \in Y$.
A part of the right hand side is calculated as
\[ \begin{array}{ll}
 \quad  \frac{1}{2 \pi i} \int_{\Gamma}~ ({\rm Log} \lambda) ~(\lambda I - U(t,s) - \kappa I)^{-2} U(t,s) ~ d \lambda  \vspace{1.5mm} \\
 = \frac{1}{2 \pi i} \int_{\Gamma} \frac{1}{\lambda}  ~(\lambda I -U(t,s)
 - \kappa I )^{-1} ~ U(t,s) ~ d \lambda \vspace{1.5mm} \\
 = \frac{1}{2 \pi i} \int_{\Gamma} \frac{1}{\lambda}  ~(\lambda I -U(t,s)- \kappa I)^{-1}  ~ \{ \lambda I - \kappa I -(\lambda I - U(t,s) - \kappa I ) \} ~ d \lambda \vspace{1.5mm} \\
 = \frac{1}{2 \pi i} \int_{\Gamma} (\lambda I -U(t,s)-\kappa I)^{-1} ~ d \lambda
 - \frac{1}{2 \pi i} \int_{\Gamma} \frac{\kappa}{\lambda} (\lambda I -U(t,s)-\kappa I )^{-1} ~
 d \lambda 
 - \frac{1}{2 \pi i} \int_{\Gamma} \frac{1}{\lambda} ~ d \lambda \vspace{1.5mm} \\
 = \frac{1}{2 \pi i} \int_{\Gamma} (\lambda I -U(t,s)-\kappa I)^{-1} ~ d \lambda
 - \frac{1}{2 \pi i} \int_{\Gamma} \frac{\kappa}{\lambda} (\lambda I -U(t,s)-\kappa I)^{-1} ~
 d \lambda \vspace{1.5mm} \\
 =  \frac{1}{2 \pi i} \int_{\Gamma} (\lambda I -U(t,s)- \kappa I )^{-1} ~ d \lambda   
 - \kappa (U(t,s)+ \kappa I )^{-1} \left\{  \frac{1}{2 \pi i} \int_{\Gamma}  \frac{1}{\lambda}  (U(t,s)+ \kappa I) (\lambda I - U(t,s)-\kappa I)^{-1} d \lambda  \right\}  \vspace{1.5mm} \\
 =    \frac{1}{2 \pi i} \int_{\Gamma} (\lambda I -U(t,s)-\kappa I )^{-1} ~ d \lambda 
 - \kappa (U(t,s)+\kappa I)^{-1} \left\{  \frac{1}{2 \pi i} \int_{\Gamma} (\lambda I - U(t,s)-\kappa I)^{-1} d \lambda -\frac{1}{2 \pi i}  \int_{\Gamma} \frac{1}{\lambda} d \lambda  \right\}  \vspace{1.5mm} \\
 =   (I - \kappa (U(t,s)+\kappa I)^{-1} ) ~  \frac{1}{2 \pi i} \int_{\Gamma} (\lambda I - U(t,s)-\kappa I)^{-1} d \lambda   \vspace{1.5mm} \\
 = (I - \kappa (U(t,s)+ \kappa I)^{-1}) ~ \frac{1}{2 \pi i} \int_{|\nu| = r} {\displaystyle \sum_{n=1}^{\infty}} \frac{U(t,s)^{n}}{\nu^{n+1}} ~ d \nu     \vspace{1.5mm} \\
  = I- \kappa (U(t,s)+ \kappa I)^{-1},  
\end{array} \] 
due to the integration by parts, where $|\lambda - \kappa| = |\nu| =r$ is a properly chosen circle large enough to include $\Gamma$.
$ (2 \pi i)^{-1}  \int_{\Gamma} \lambda^{-1} d \lambda = 0$ is seen by applying $d {\rm Log} \lambda/d \lambda = 1/\lambda$.  
$(2 \pi i)^{-1} \int_{|\nu| = r} {\displaystyle \sum_{n=1}^{\infty}} U(t,s)^{n} \nu^{-n-1} ~ d \nu  = I$ follows from the singularity of $\nu^{-n-1}$. 

Consequently
\begin{equation} \label{intmed} \begin{array}{ll}
A(t) ~ u_s = \{I- \kappa (U(t,s)+\kappa I)^{-1}\}^{-1} ~ \partial_{t} {\rm Log} ~ (U(t,s) + \kappa I) ~ u_s  \vspace{1.5mm} \\
\quad = (U(t,s)+\kappa I) U(t,s)^{-1} ~ \partial_{t} {\rm Log} ~ (U(t,s) + \kappa I) ~ u_s  \vspace{1.5mm} \\
\quad =  (I+\kappa U(s,t))~ \partial_{t} {\rm Log} ~ (U(t,s) + \kappa I) ~ u_s 
\end{array} \end{equation}
is obtained for $u_s \in Y$. \quad $\square$ \\


The meaning of logarithmic representation is examined by focusing on $\partial_{t} {\rm Log} ~ (U(t,s) + \kappa I)$.
What is introduced by Eq.~\eqref{logex} is a kind of resolvent approximation\index{resolvent approximation} of $A(t)$
\[ \begin{array}{ll}
 \partial_{t} {\rm Log} ~ (U(t,s) + \kappa I) =   (I+\kappa U(s,t))^{-1} A(t), 
\end{array} \]
in which $A(t)$ is represented by the resolvent operator of $U(s,t)$.
As seen in the following it is notable that there is no need to take $\kappa \to 0$.
This point is different from the usual treatment of resolvent approximations; indeed it is impossible to take  $\kappa \to 0$ if the origin is not included in the resolvent set of $U(t,s)$.
On the other hand, it is also seen by Eq.~\eqref{logex} that
\begin{equation}   \label{simtra}
 \partial_{t} {\rm Log} ~ (U(t,s) + \kappa I) =  (U(t,s)+\kappa I )|_{\kappa = 0} A(t)  (U(t,s)+\kappa I)^{-1}     
\end{equation}
shows a structure of similarity transform\index{similarity transform}, where $(U(t,s)+\kappa)|_{\kappa=0}$ means $U(t,s)+\kappa$ satisfying a condition $\kappa=0$.
This asymmetric similarity transform from left and right hand sides are remarkable, and it becomes symmetric if $\kappa=0]$.
A part $\partial_{t} {\rm Log} ~ (U(t,s) + \kappa I)$ plays an essential role in the following discussion.

\section{Regularized evolution operator}
\subsection{Alternative infinitesimal generator}
The alternative infinitesimal generator\index{alternative infinitesimal generator} is introduced in order to extract bounded parts from the pre-infinitesimal generator $A$ \cite{17iwata-3}. 
The operator $A(t) \in G(X)$ is generally unbounded in $X$.
A bounded operator $a(t,s)$ on $X$ is introduced.  \vspace{2.5mm}\\

\begin{definition}[Alternative infinitesimal generator]
Let $\kappa$ be a certain complex number.
For a certain $v_s \in X$, the alternative infinitesimal generator $\partial_{t} a(t,s)$ to $A(t) \in G(X)$ is defined using
\begin{equation} \label{logex3} \begin{array}{ll}
a(t,s) v_s := {\rm Log} ~ (U(t,s) + \kappa I) v_s
\end{array} \end{equation}
on $X$, where $\partial_t $ denotes  $t-$differential in a weak sense.
\end{definition}

In the present setting assuming the existence of $\partial_t U(t,s)$ and therefore $A(t)$, the operator $\partial_t a(t,s)$ exists.
According to the logarithmic representation,
\[
 A(t) ~ u_s =  (I - \kappa (U(t,s) + \kappa I)^{-1})^{-1} ~ \partial_{t} a(t,s) ~ u_s
\]
is obtained.
Since $a(t,s)$ is a bounded operator defined by the Dunford-Riesz integral, $v_s$ in the definition of alternative infinitesimal generator can be taken from $X$.
Since $\kappa$ is chosen to separate the spectral set of $U(t,s)+\kappa I$ from the origin, the inverse operator of $U(t,s) + \kappa I$ always exists, as $\kappa$ is taken from $\{ \lambda \in {\mathbb C}; ~ |\lambda| > M e^{\beta T} \}$.

\begin{definition}[Regularized evolution operator]
The alternative infinitesimal generator $\partial_{t} a(t,s)$ generates the regularized evolution operator
\begin{equation}  \begin{array}{ll}
 e^{a(t,s)}  = {\displaystyle \sum_{n=0}^{\infty}} \frac{ a(t,s)^n}{n!},
\end{array} \end{equation}
which is represented by the convergent power series.
\end{definition}

The operator $ e^{a(t,s)} $ is regularized in the following sense; the inverse evolution operator $e^{-a(t,s)} $ always exists, if $ e^{a(t,s)} $ exists.
This fact, which arises from the boundedness of $a(t,s)$, is true, even if negative time evolution $U(t,s)^{-1}= U(s,t)$ is not well-defined and only positive time evolution $U(t,s)$ is given ($t>s$).

It is remarkable that
\[ 
{\rm Log}  e^{a(t,s)}  = {\rm Log} (U(t,s)+\kappa I) = a(t,s)
\]
is always satisfied, while 
\[ \int  (I+ \kappa U(s,t))^{-1} A(t) ~ dt  =  {\rm Log} ~  (U(t,s) + \kappa I)  \]
is not necessarily satisfied because the limited range of imaginary spectral distribution is necessarily true only for the right hand side.
In this sense $a(t,s)$ corresponds to the extracted bounded part of the infinitesimal generator $A(t)$.
The regularized trajectory\index{regularized trajectory} in finite/infinite dimensional dynamical systems\index{dynamical system} (for textbooks, see \cite{74hirsch,97temam}) arises from the regularized evolution operator.
Note that, as the well-defined $\partial_t a(t,s)$ is not necessary for $e^{a(t,s)}$ to be well-defined, only the well-defined $a(t,s)$ is sufficient for $e^{a(t,s)}$ to be well-defined.
This fact essentially simplifies the discussion in applying $e^{a(t,s)}$ in which there is no need to consider weak differential. 

Using the relation between the logarithm and the exponential functions,
\[  U(t,s) =  e^{a(t,s)}  - \kappa I \]
is valid.
It shows a correspondence between $ e^{a(t,s)} $ and $U(t,s)$ at the level of an evolution operator.
One difference is whether the semigroup property is satisfied or not, and another difference is whether the convergence power series representation is always true or not.
Meanwhile, at the level of infinitesimal generators, there is a substantial difference between $a(t,s) $ and $A(t)$. 
That is, $a(t,s)$ is always bounded on $X$, while $A(t)$ is not necessarily bounded on $X$.
Because of the boundedness of $a(t,s)$ on $X$, the inverse operator $e^{-a(t,s)}$ always exists if $e^{a(t,s)}$ exists.
One of the essential ideas is to generate $ e^{a(t,s)}$, instead of generating $U(t,s)$.

\begin{theorem}[Modified semigroup property]
Let $\kappa$ be a certain complex number.
For the operator $e^{a(t,s)}$ on $X$, the semigroup property is replaced with
\begin{equation} \begin{array}{ll} 
 e^{a(t,s)} -  e^{a(t,r)} e^{a(r,s)}  =  \kappa (\kappa  + 1) I  - \kappa  (e^{a(t,r)} + e^{a(r,s)} ), \vspace{1,5mm}  \\
 e^{a(s,s)}  - I    =   \kappa  I.
 \label{theq01} \end{array}  \end{equation}
The inverse relation is replaced with
\begin{equation} \begin{array}{ll}
  e^{a(s,t)} e^{a(t,s)}  - e^{a(s,s)}  =  \kappa  (e^{a(t,s)} + e^{a(s,t)} ) -  \kappa (\kappa + 1) I . 
 \label{theq02} \end{array} \end{equation}
In particular the commutation
\begin{equation} \begin{array}{ll}
  e^{a(s,t)} e^{a(t,s)}  - e^{a(t,s)}  e^{a(s,t)}   = 0
 \label{theq03} \end{array} \end{equation}
is necessarily valid.
\label{thm:semi}
\end{theorem}

{\bf Proof. ~}
Substitution of $U(t,s) = e^{a(t,s)} - \kappa I$ to $U(t,r)~U(r,s) = U(t,s)$ leads to the following relation:
\[ \begin{array}{ll}
U(t,r)~U(r,s) =  e^{a(t,s)} - \kappa I,   
\end{array} \]
and 
\[ \begin{array}{ll}
( e^{a(t,r)} - \kappa I)~( e^{a(r,s)} - \kappa I) =  e^{a(t,s)} - \kappa I.   
\end{array} \]
where, by taking $\kappa$ with a large $|\kappa|$, $\kappa$ is possible to be taken as common to $U(t,s)$ with different $t$ and $s$.
Meanwhile the replacement of $U(t,s) = e^{a(t,s)} - \kappa I$ with $U(s,s) = I$ leads to the following relation:
\[ \begin{array}{ll}
 e^{a(s,s)}  =  ( \kappa  +1) I.  
\end{array} \]
That is, for $\kappa \ne 1$, $( \kappa  +1)^{-1} e^{a(s,s)}$ behaves as the unit operator.
Modified version of semigroup property (i.e., \eqref{theq01}) has been proved.
The inverse relation \eqref{theq02} follows readily from Eq.~\eqref{theq01}.
According to Eq.~\eqref{theq01},
\[ \begin{array}{ll} 
 e^{a(t,t)} -  e^{a(t,s)} e^{a(s,t)}  =  \kappa (\kappa  + 1) I  - \kappa  (e^{a(t,s)} + e^{a(s,t)} ), \vspace{1,5mm}  \\
\end{array}  \]
is valid. 
Combination with another relation
\[ \begin{array}{ll} 
 e^{a(s,s)} -  e^{a(s,t)} e^{a(t,s)}  =  \kappa (\kappa  + 1) I  - \kappa  (e^{a(s,t)} + e^{a(t,s)} ), \vspace{1,5mm}  \\
\end{array}  \]
leads to the commutation:
\[ \begin{array}{ll} 
   e^{a(t,s)} e^{a(s,t)}
 - e^{a(s,t)} e^{a(t,s)} 
 = 0
\end{array}  \]
where $ e^{a(t,t)} = e^{a(s,s)} = (1+\kappa)I$ is utilized.  \quad $\square$ \\

Equations~\eqref{theq01} and \eqref{theq02} show the commutativity and violation of semigroup property by $e^{a(t,s)}$.
The right hand sides of Eqs.~\eqref{theq01} and \eqref{theq02} are equal to zero for $\kappa = 0$.
These situations correspond to the cases when the semigroup property is satisfied by $e^{a(t,s)}$, and it is readily seen that the insufficiency of semigroup property arises from the introduction of nonzero $\kappa$.

The decomposition is obtained by the following structure theorem for the regularized evolution operator.
Note that the decomposition of $e^{a(t,s)}$ also provides a certain relation between the time-discretization and the violation of semigroup property.

\begin{theorem}[Structure of regularized evolution operator]
Let $\kappa$ be a certain complex number.
For a given decomposition $s< r_1, r_2, \cdots , r_n <t$ of the interval $[s,t]$ with $n \ge 2$, the operator $e^{a(t,s)}$ on $X$ is represented by
\begin{equation} \begin{array}{ll}
 e^{a(t,s)} 
    =    e^{a(t,r_n)} e^{a(r_{n},r_{n-1})} \cdots  e^{a(r_2,r_1)}  e^{a(r_1,s)}  \vspace{1.5mm} \\
\quad  +  \kappa (\kappa  + 1) I  
    - \kappa   (e^{a(t ,r_{1})} + e^{a(r_{1},s)} ) \\
\quad +  {\displaystyle \sum_{k=2}^n}  \left[   \{  \kappa (\kappa  + 1)    
    -    \kappa  (e^{a(t ,r_{k})} + e^{a(r_{k},r_{k-1})} ) \}
     ~ e^{a(r_{k-1},r_{k-2})} \cdots   e^{a(r_2,r_1)}  e^{a(r_1,s)}    \right]  
 \label{theq04} \end{array} \end{equation}
where $r_0$ and $r_{n+1}$ in the sum are denoted as $s = r_0$ and $t = r_{n+1}$ respectively.
\end{theorem}

{\bf Proof. ~}
According to Eq.~\eqref{theq01}, a decomposition
\[ \begin{array}{ll}
 e^{a(t,s)}  =  e^{a(t,r_1)} e^{a(r_1,s)}  + \kappa (\kappa  + 1) I  - \kappa  (e^{a(t,r_1)} + e^{a(r_1,s)} )
\end{array} \]
is true.
Another decomposition
\[ \begin{array}{ll}
 e^{a(t,r_1)} =  e^{a(t,r_2)} e^{a(r_2,r_1)} +  \kappa (\kappa  + 1) I  - \kappa  (e^{a(t,r_2)} + e^{a(r_2,r_1)} ),
\end{array} \]
is also true, and then
\[ \begin{array}{ll}
 e^{a(t,s)}
  =    \{ e^{a(t,r_2)} e^{a(r_2,r_1)} +  \kappa (\kappa  + 1) I  - \kappa  (e^{a(t,r_2)} + e^{a(r_2,r_1)} ) \}
  e^{a(r_1,s)} +  \kappa (\kappa  + 1) I  - \kappa  (e^{a(t,r_1)} + e^{a(r_1,s)} )  \vspace{2.5mm} \\
  =    e^{a(t,r_2)} e^{a(r_2,r_1)}  e^{a(r_1,s)}     +  \kappa (\kappa  + 1) (I  +  e^{a(r_1,s)} )    
  - \kappa  \{ e^{a(t,r_1)} + e^{a(r_1,s)}  +  ( e^{a(t,r_2)} + e^{a(r_2,r_1)} )   e^{a(r_1,s)} \}  
\end{array} \]
follows by sorting based on $\kappa$ and $ \kappa (\kappa  + 1)$ dependence.
Further decomposition shows
\[ \begin{array}{ll}
 e^{a(t,r_2)} =  e^{a(t,r_3)} e^{a(r_3,r_2)} +  \kappa (\kappa  + 1) I  - \kappa  (e^{a(t,r_3)} + e^{a(r_3,r_2)} ),
\end{array} \]
and then
\[ \begin{array}{ll}
 e^{a(t,s)} 
  =    e^{a(t,r_2)} e^{a(r_2,r_1)}  e^{a(r_1,s)}   \\
  +  \kappa (\kappa  + 1) (I  +  e^{a(r_1,s)} )    - \kappa  \{ e^{a(t,r_1)} + e^{a(r_1,s)}  +  ( e^{a(t,r_2)} + e^{a(r_2,r_1)} )   e^{a(r_1,s)} \} \vspace{2.5mm} \\
  =    \{   e^{a(t,r_3)} e^{a(r_3,r_2)} +  \kappa (\kappa  + 1) I  - \kappa  (e^{a(t,r_3)} + e^{a(r_3,r_2)} )  \} e^{a(r_2,r_1)}  e^{a(r_1,s)}   \\
  +  \kappa (\kappa  + 1) (I  +  e^{a(r_1,s)} )    - \kappa  \{ e^{a(t,r_1)} + e^{a(r_1,s)}  +  ( e^{a(t,r_2)} + e^{a(r_2,r_1)} )   e^{a(r_1,s)} \} \vspace{2.5mm} \\
    =    e^{a(t,r_3)} e^{a(r_3,r_2)}  e^{a(r_2,r_1)}  e^{a(r_1,s)}  +  \kappa (\kappa  + 1)  e^{a(r_2,r_1)}  e^{a(r_1,s)}  - \kappa  (e^{a(t,r_3)} + e^{a(r_3,r_2)} )    e^{a(r_2,r_1)}  e^{a(r_1,s)}   \\
  +  \kappa (\kappa  + 1) (I  +  e^{a(r_1,s)} )    - \kappa  \{ e^{a(t,r_1)} + e^{a(r_1,s)}  +  ( e^{a(t,r_2)} + e^{a(r_2,r_1)} )   e^{a(r_1,s)} \}    \vspace{2.5mm} \\
    =    e^{a(t,r_3)} e^{a(r_3,r_2)}  e^{a(r_2,r_1)}  e^{a(r_1,s)}   
  +  \kappa (\kappa  + 1) \{ I  +  (I +  e^{a(r_2,r_1)})  e^{a(r_1,s)}  \}   \\
   - \kappa  \{ e^{a(t,r_1)} + e^{a(r_1,s)}  +  ( e^{a(t,r_2)} + e^{a(r_2,r_1)} )   e^{a(r_1,s)}
   +  (e^{a(t,r_3)} + e^{a(r_3,r_2)} )    e^{a(r_2,r_1)}  e^{a(r_1,s)}   \}  
\end{array} \]
follows.
For a certain $n \ge 2$, a constitutional representation is suggested by the deduction:
\[ \begin{array}{ll}
 e^{a(t,s)} 
    =    e^{a(t,r_n)} e^{a(r_{n},r_{n-1})} \cdots  e^{a(r_2,r_1)}  e^{a(r_1,s)}   \\
  +  \kappa (\kappa  + 1)  \left[ I + e^{a(r_1,s)} + \left( e^{a(r_2,r_1)} e^{a(r_1,s)} \right)  + \cdots  + \left( e^{a(r_n,r_{n-1})} \cdots e^{a(r_2,r_1)} e^{a(r_1,s)}  \right) \right]   \\
   - \kappa  \left[ e^{a(t,r_1)} + e^{a(r_1,s)} 
   +  ( e^{a(t,r_2)} + e^{a(r_2,r_1)} )   e^{a(r_1,s)}  \right. \\
   \left. +  (e^{a(t,r_3)} + e^{a(r_3,r_2)} )    e^{a(r_2,r_1)}  e^{a(r_1,s)}  \cdots
    +  (e^{a(t,r_n)} + e^{a(r_{n},r_{n-1})} )  e^{a(r_{n-1},r_{n-2})} \cdots   e^{a(r_2,r_1)}  e^{a(r_1,s)} 
     \right].   
\end{array} \]
Consequently
\[ \begin{array}{ll}
 e^{a(t,s)} 
    =    e^{a(t,r_n)} e^{a(r_{n},r_{n-1})} \cdots  e^{a(r_2,r_1)}  e^{a(r_1,s)}   \\
\quad  +  \kappa (\kappa  + 1) \left( I + {\displaystyle \sum_{k=2}^n}  \left[ e^{a(r_{k-1},r_{k-2})} \cdots e^{a(r_2,r_1)} e^{a(r_1,s)}  \right] \right)   \\
\quad    - \kappa
\left(   (e^{a(t ,r_{1})} + e^{a(r_{1},s)} )
	  +   {\displaystyle \sum_{k=2}^{n}}  \left[  (e^{a(t ,r_{k})} + e^{a(r_{k},r_{k-1})} ) ~ e^{a(r_{k-1},r_{k-2})} \cdots   e^{a(r_2,r_1)}  e^{a(r_1,s)}    \right]  \right)
\end{array} \]
is obtained.
The statement is proved by sorting terms. \quad $\square$ \\

Using the regularized evolution operator, the logarithmic representation is readily generalized to the infinitesimal generators of invertible and  non-invertible evolution operators.
Indeed, according to the proof of Theorem \ref{thm1}, only the boundedness of $U(t,s)$ on $X$ and the resulting time-interval symmetry is essential.

\begin{corollary} [Generalized logarithmic representation of infinitesimal generators]
\label{cor1}
Let $t$ and $s$ satisfy $-T \le t,s \le T$, and $Y$ be a dense subspace of $X$.
For non-invertible $U(t,s)$; $U(t,s)$ defined in Sec.~\ref{tp-group} without assuming
\[ U(t,s)^{-1} = U(s,t), \]
 let $A(t) \in G(X)$ and $\partial_t U(t,s)$ be determined by Eqs.~\eqref{pe-group} and \eqref{de-group} respectively.  
If $A(t)$ and $U(t,s)$ commute, pre-infinitesimal generators $\{ A(t) \}_{-T \le t \le T}$ are represented by means of the logarithm function; there exists a certain complex number $\kappa \ne 0$ such that
\begin{equation} \begin{array}{ll}
 A(t) ~ u_s =  (I - \kappa e^{-a(t,s)})^{-1} ~ \partial_{t} {\rm Log} ~ (U(t,s) + \kappa I) ~ u_s, 
\end{array} \end{equation}
where $u_s$ is an element in $Y$.
\end{corollary}

{\bf Proof. ~ } 
The first line of Eq.~(\ref{intmed}) shows the validity of the statement. 
Indeed, for the logarithmic representation, the invertible property does not play any roles after introducing nonzero $\kappa \in {\mathbb C}$.
In particular $  (I - \kappa e^{-a(t,s)})^{-1} $ is always well defined for a certain $\kappa$.
 \quad $\square$ \\

Uisng the regularized evolution operator,
A similarity transform representation~(\ref{simtra}) for $\partial_{t}  a(t,s)$ is written by
\begin{equation} \label{simtra2}
 \partial_{t} a(t,s) =  e^{a(t,s)}|_{\kappa = 0} A(t)  e^{-a(t,s)},    
\end{equation}
where the boundedness of $a(t,s)$ allows us to define $e^{-a(t,s)}$.
Due to the boundedness of $a(t,s)$ on $X$, $ e^{a(t,s)} $ is always well-defined by a convergent power series\index{convergent power series}.
It leads to the holomorphic property of $ e^{a(t,s)} $.
Here is the reason why $ e^{a(t,s)}$ is called the regularized evolution operator\index{regularized evolution operator}. 
 
In the following, the generalized  logarithmic representation of infinitesimal generators utilized.
It enables us to have the logarithmic representation not only for the $C_0$-groups but for the $C_0$-semigroups.

\subsection{Renormalized abstract evolution equations} 
Evolution equations are renormalized by means of the alternative infinitesimal generators and regularized evolution operator.
\begin{corollary}[Renormalized abstract evolution equations]  \label{transform}
If $a(t,s)$ with different $t$ and $s$ are further assumed to commute,
\begin{equation} \label{replce} \begin{array}{ll}
\partial_t  e^{a(t,s)} v_s = [ \partial_t a(t,s) ] ~   e^{a(t,s)} v_s 
\quad {\rm leading~to} \quad  \partial_t v(t) = [ \partial_t a(t,s) ]  v (t)  \\
\end{array} \end{equation}
is satisfied for $v_s \in Y$, where $\partial_t$ denotes $t-$differential in a weak sense.
This is a linear evolution equation satisfied by $v(t) = e^{a(t,s)}  v_s$.
\end{corollary}

{\bf Proof. ~}
For the evolution operator
\begin{equation} \label{convp} \begin{array}{ll} 
e^{a(t,s)} = {\displaystyle \sum_{n=0}^{\infty}} \frac{ a(t,s)^n}{n!}, \end{array} \end{equation}
the existence of $\partial_t e^{a(t,s)}$ is ensured by the existence of $\partial_t a(t,s)$.
Using the commutation between $a(t,s)$ with different $t$ and $s$, 
\[ \begin{array}{ll}
\partial_t \{ a(t,s) \}^n 
 = n \{ a(t,s) \}^{n-1} ~ \partial_t a(t,s)
\end{array} \] 
is true, and the homogeneous-type abstract evolution equation is rephrased as an equation with bounded infinitesimal generator
\begin{equation} \label{traeq} \begin{array}{ll}
\partial_t e^{a(t,s)} v_s
 = e^{a(t,s)} ~ (\partial_t a(t,s))  v_s 
\end{array} \end{equation}
for $v_s \in Y$. \quad $\square$ \\

This is an abstract evolution equation obtained by the replacement $A(t)$ with $\partial_t a(t,s)$.
Here one essential idea is to generate $ e^{a(t,s)} $ instead of $U(t,s)$; although $e^{a(t,s)}$ is easily defined due to the boundedness of $a(t,s)$, the general unboundedness of infinitesimal generator $\partial_t a(t,s)$ in $X$ is ensured by the similarity transform~(\ref{simtra2}).
Under the commutation assumption between $A(t)$ and $U(t,s)$, Eq.~\eqref{traeq} is rephrased as
\[ \begin{array}{ll}
\partial_t e^{a(t,s)} v_s
 =   \partial_t (U(t,s)+\kappa I) v_s
 =   \partial_t U(t,s) v_s,
\end{array} \]
and 
\[ \begin{array}{ll}
  e^{a(t,s)} ~ (\partial_t a(t,s) ) v_s 
 = (U(t,s)+\kappa I) \partial_t( {\rm Log} (U(t,s)+\kappa I)) v_s \vspace{1.5mm} \\ 
  = (U(t,s)+\kappa I)  (U(t,s)+\kappa I)^{-1} U(t,s) A(t) v_s \vspace{1.5mm} \\
 =  U(t,s) A(t) v_s  \vspace{1.5mm} \\
=  A(t) U(t,s) v_s,
\end{array} \]
where Theorem~\ref{thm1} is applied.
Consequently
\[  \partial_t U(t,s) v_s =   A(t) U(t,s) v_s \]
is obtained.
Note again that $e^{a(t,s)}$ does not satisfy the semigroup property, while $U(t,s)$ satisfies it.

\subsection{Linearized infinitesimal generator}
The linearity of the semigroup is not assumed in the preceding discussion, so that the operator $U(t,s)$ can be taken as either linear or nonlinear semigroup.
Let us assume a more general situation, in which 
\begin{itemize}
\item the existence of $U(t,s)$ satisfying Eq.~(\ref{sg1}) is locally ture for $t,s \in [-T, T]$,  
\item  the existence of the infinitesimal generator of $U(t,s)$ is not clear.
\end{itemize}
This situation corresponds to the situation when only the unique local-existence of their solutions the nonlinear partial differential equations is ensured.
One of the important application of the renormalized abstract evolution equation\index{renormalized abstract evolution equation} is the linearization, which enables to analyze the local profile of nonlinear semigroup.

\begin{corollary}[Linearized evolution equation]
For $-T \le t, s  \le T$, the two-parameter group is defined on $X$.
Let either linear or nonlinear semigroup $U(t,s)$ defined on a Banach space $X$ satisfy Eq.~(\ref{sg1}).
Let $U(t,s) u_s$ be the solution of nonlinear equation: 
\[ \partial_t u(t,s) = F(u(t,s)) .   \]    
If  the logarithmic representation ${\rm Log}(U(t,s) + \kappa I)$ is true, 
\[ \begin{array}{ll} 
 \partial_t v(t) = [ (I - \kappa (U(t,s) + \kappa I)^{-1})^{-1}  \partial_t a(t,s) ] ~ v(t)  
\end{array} \]
is the linearized equation, where note that $a(t,s)$ includes a parameter $\kappa \in {\mathbb C}$.
If the logarithmic representation is true for $\kappa =0$, the infinitesimal generator of linearized problem is simply represented by $\partial_t a(t,s)|_{\kappa = 0} = \partial_t {\rm Log}(U(t,s))$. 
\end{corollary} 

The condition for obtaining the linearized evolution equation\index{linearized evolution equation} is the locality for the evolution direction $t$, which leads to the boundedness of the spectral set of $U(t,s)$. 
The theoretical procedure of obtaining the linearized problem is summarized as follows.
For nonzero $\kappa \in {\mathbb C}$, first, $U(t,s)$ is regarded as an exponential function; second, calculating  the logarithm of $U(t,s)$; and finally the linearized operator
\[ \begin{array}{ll} 
\quad (I - \kappa (U(t,s) + \kappa I)^{-1})^{-1} ~ \partial_{t} a(t,s)  \vspace{1.5mm} \\
= \partial_{t} [ (I - \kappa (U(t,s) + \kappa I)^{-1})^{-1} ~ a(t,s) ]
-  \partial_{t} [(I - \kappa (U(t,s) + \kappa I)^{-1})^{-1} ] ~ a(t,s)
\end{array} \]
generates the regularized evolution operator
\[ \begin{array}{ll} 
  e^{a(t,s)} - \kappa I.
\end{array} \]
It is more clearly understood by the case of $\kappa = 0$,
\[ \begin{array}{ll} 
``U(t,s) = e^{a(t,s)}|_{\kappa =0} "
\quad {\rm leads~to} \quad
`` {\rm Log}e^{a(t,s)}|_{\kappa =0} 
= a(t,s)|_{\kappa =0} ".
\end{array} \]
Consequently the operator-logarithm is regarded as a mapping from ``continuous group'' to ``bounded algebra''.
These alternative equations can be used to analyze quasi-linear evolution equations\index{quasi-linear evolution equations} and full-nonlinear evolution equations\index{full-nonlinear evolution equations}~\cite{75kato}.
\vspace{3mm} \\

\subsubsection{Autonomous case} \label{homosection}
The regularity results have not been much studied in the Cauchy problem of hyperbolic partial differential equations (for a textbook, see \cite{85mizohata}).
The regularized evolution operator, which is also applicable to some hyperbolic type equations, is utilized to solve autonomous Cauchy problems.  
\begin{equation} \label{homoporo} \left\{  \begin{array}{ll}
\partial_t u(t)  = A(t) u(t) \vspace{2.5mm} \\
u(s) = u_s, 
\end{array} \right. \end{equation} 
in $X$, where $A(t) \in G(X):Y \to X$ is assumed to be an infinitesimal generator of $U(t,s)$ satisfying the semigroup property, $-T \le t,s \le T$ is satisfied, $Y$ is a dense subspace of $X$ permitting the representation shown in Eq.~\eqref{logex}, and $u_s$ is an element of $X$.

As seen in Eq.~\eqref{replce}, under the assumption of commutation, a related Cauchy problem is obtained as
\begin{equation} \left\{  \begin{array}{ll} \label{reweq}
\partial_t v(t,s)  = (\partial_t a(t,s)) ~ v(t,s) \vspace{2.5mm} \\
v(s,s) = e^{a(s,s)} u_s,
\end{array} \right. \end{equation}
in $X$, where $\partial_t a(t,s) = \partial_t {\rm Log} (U(t,s)+\kappa I)$ is well-defined.
It is possible to solve the re-written Cauchy problem, and the solution is represented by
\[  \begin{array}{ll}
v(t,s) =  e^{a(t,s)} u_s  =
 {\displaystyle \sum_{n=0}^{\infty} } \frac{ a(t,s)^n}{n!} u_s
 \end{array} \]
 for $u_s \in X$  (cf.~Eq.~\eqref{convp}).

\begin{theorem} \label{hols}
Operator $e^{a(t,s)}$ is holomorphic.
\end{theorem}

{\bf Proof ~ }
According to the boundedness of $a(t,s)$ on $X$ (Lemma~\ref{lem3}), $\partial_t^n  e^{a(t,s)}$~\cite{51taylor} is possible to be represented as
\begin{equation} \label{anreap} \begin{array}{ll}
\partial_t^n  e^{a(t,s)} = \frac{1}{2 \pi i} \int_{\Gamma} \lambda^n e^{\lambda} (\lambda I -a(t,s))^{-1} ~ d \lambda, 
\end{array} \end{equation} 
for a certain $\kappa \in {\mathbb C}$, where $ \lambda^n e^{\lambda}$ does not hold any singularity for any finite $\lambda$.
Following the standard theory of evolution equation, 
\[  \begin{array}{ll}
 \| \partial_t^n  e^{a(t,s)}   \| \le \frac{C_{\theta,n}}{ \pi (t \sin \theta)^n}
 \end{array} \] 
is true for a certain constant $C_{\theta,n}$  ($n = 0,1,2,\cdots$), where $\theta \in (0 \pi/2)$ and $|\arg t| < \pi/2$ are satisfied (for the detail, e.g., see \cite{79tanabe}).
It follows that
\begin{equation} \label{leiq} \begin{array}{ll}
{ \displaystyle  \lim_{t \to +0} } \sup t^n \| \partial_t^n  e^{a(t,s)} \|
\le {\displaystyle \lim_{t \to +0}} \sup t^n   \frac{C_{\theta,n}}{\pi (t \sin \theta)^n}  < \infty.
\end{array} \end{equation} 
Consequently, for $|z-t|<t \sin \theta$, the power series expansion
\[ \begin{array}{ll}
{\displaystyle \sum_{n=0}^{\infty} } \frac{(z-t)^n}{n !}   \partial_t^n e^{a(t,s)} 
 \end{array} \]
is uniformly convergent in a wider sense.
Therefore $e^{a(t,s)}$ is holomorphic. \quad $\square$ \\

\begin{theorem} \label{reprr}
For $u_s \in X$ there exists a unique solution $u(\cdot) \in C([-T,T];X)$ for \eqref{homoporo} with a convergent power series representation:
\begin{equation} \label{dairep} \begin{array}{ll} 
 u(t) = U(t,s) u_s =  ( e^{a(t,s)}  - \kappa  I ) u_s = \left( {\displaystyle \sum_{n=0}^{\infty} } \frac{ a(t,s)^n}{n!}  - \kappa I \right)  u_s,
\end{array} \end{equation}
where $\kappa$ is a certain complex number.
\end{theorem}

{\bf Proof ~}
The unique existence follows from the assumption for $A(t)$.
The regulariized evolution operator $e^{a(t,s)}$ is holomorphic function (Theorem~\ref{hols}) with the convergent power series representation (Eq.~\eqref{convp}).
By applying $U(t,s) + \kappa I  = e^{a(t,s)} $ the solution of the original Cauchy problem is obtained as
 \[ \begin{array}{ll}
  u(t) =   ( e^{a(t,s)}  - \kappa  I ) u_s = \left( {\displaystyle \sum_{n=0}^{\infty} } \frac{ a(t,s)^n}{n!}  - \kappa I \right)  u_s,
  \end{array} \]
for the initial value $u_s \in X$.
Note that $A(t)$ is not assumed to be a generator of analytic evolution family, but only a generator of evolution family \quad \quad $\square$ \\

For $I_{\lambda}$ denoting the resolvent operator of $A(t)$, the evolution operator defined by the Yosida approximation\index{Yosida approximation} is written by
\[ \begin{array}{ll}
u(t) = {\displaystyle \lim_{\lambda \to 0} } \exp ( \int_s^t I_{\lambda} A(\tau) ~d \tau ) u_s,
\end{array} \] 
so that more informative representation is provided by Theorem~\ref{reprr} compared to the standard theory based on the Hille-Yosida theorem\index{Hille-Yosida theorem}.

\subsubsection{Non-autonomous case}
Series representation in autonomous part leads to the enhancement of the solvability.
Let $Y$ be a dense subspace of $X$ permitting the representation shown in Eq.~\eqref{logex},  and $u_s$ is an element of $X$.
The regularized evolution operator is utilized to solved non-autonomous Cauchy problems.
\begin{equation} \label{origih} \left\{  \begin{array}{ll}
\partial_t u(t)  = A(t)  u(t) + f(t) \vspace{2.5mm} \\
u(s) = u_s
\end{array} \right. \end{equation}
in $X$, where $A(t) \in G(X):Y \to X$ is assumed to be an infinitesimal generator of $U(t,s)$ satisfying the semigroup property, $f \in L^1(-T,T;X)$ is locally H\"older continuous on $[-T,T]$
\[ \begin{array}{ll} 
\| f(t) - f(s) \| \le C_{H} |t-s|^{\gamma}
\end{array} \]
for a certain positive constant $C_{H}$, $\gamma \le 1$ and $-T \le t,s \le T$. 
The solution of non-autonomous problem does not necessarily exist in such a setting (in general, $f \in C([-T,T];X)$ is necessary).

\begin{theorem} \label{thm-inh}
Let $f \in L^1(-T,T;X)$ be locally H\"older continuous on $[-T,T]$.
For $u_s \in X$ there exists a unique solution $u(\cdot) \in C([-T,T];X)$ for \eqref{origih} such that
\[ \begin{array}{ll} 
 u(t) = \left[  {\displaystyle \sum_{n=0}^{\infty}} \frac{ a(t,s)^n}{n!}  - \kappa I \right]  u_s +    {\displaystyle \int_s^t}  \left[ {\displaystyle \sum_{n=0}^{\infty} } \frac{ a(t,\tau)^n}{n!}  - \kappa I \right]  f(\tau) d\tau
\end{array} \]
using a certain complex number $\kappa$.
\end{theorem}

{\bf Proof ~}
Let us begin with cases with $f \in C([-T,T];X)$.
The unique existence follows from the standard theory of evolution equation.
The representation follows from that of $U(t,s)$ and the Duhamel's principle
\begin{equation}  \begin{array}{ll}
 u(t)  =   U(t,s) u_s + \int_s^tU(t,\tau) f(\tau) d\tau \vspace{1.5mm}\\
 =   ( e^{a(t,s)} - \kappa I )  u_s + \int_s^t  [  e^{a(t,\tau)} - \kappa I ]  f(\tau) d\tau,
\end{array}  \end{equation}
where the convergent power series representation of $e^{a(t,s)}$ is valid (cf.~Eq.~\eqref{convp}). 
\begin{equation}  \begin{array}{ll}
 u(t)   =   ( e^{a(t,s)} - \kappa I )  u_s + \int_s^t  [  e^{a(t,\tau)} - \kappa I ]  f(\tau) d\tau, \\
  =    e^{a(t,s)}   u_s + \int_s^t  [  e^{a(t,\tau)}  ]  f(\tau) d\tau
    - \kappa I  u_s + \int_s^t  [  - \kappa I ]  f(\tau) d\tau. \\
\end{array}  \end{equation}

Next let us consider cases with the locally H\"older continuous $f(t)$.
According to the linearity of Eq.~\eqref{origih}, it is sufficient to consider the inhomogeneous term.
For $\epsilon$ satisfying $0 < \epsilon << T$,
\[   \begin{array}{ll}
  \int_s^{t+\epsilon} [e^{a(t,\tau)} -\kappa I] f(\tau) d\tau  
  ~\to~  
   \int_s^{t} [e^{a(t,\tau)} -\kappa I] f(\tau) d\tau
\end{array}  \]
is true by taking $\epsilon \to 0$.
On the other hand, 
\begin{equation} \label{conveq}  \begin{array}{ll}
   A(t)  \int_s^{t+\epsilon} [e^{a(t,\tau)} -\kappa I] f(\tau) d\tau 
=  \int_s^{t+\epsilon} A(t) U(t,\tau) f(\tau) d\tau  \vspace{1.5mm}\\
=   \int_s^{t+\epsilon} A(t)  U(t,\tau)( f(\tau) - f(t)) d\tau  + \int_s^{t+\epsilon} A(t)  U(t,\tau) f(t) d\tau  \vspace{1.5mm}\\
=   \int_s^{t+\epsilon} A(t)  U(t,\tau) ( f(\tau) - f(t)) d\tau  - \int_s^{t+\epsilon}  \partial_{\tau} U(t,\tau) f(t) d\tau   \vspace{1.5mm}\\
=   \int_s^{t+\epsilon} A(t)  U(t,\tau) ( f(\tau) - f(t)) d\tau  -  U(t,t+\epsilon) f(t) + U(t,s) f(t)  \vspace{1.5mm}\\
=   \int_s^{t+\epsilon}  (I+\kappa U(s,t))\partial_{t} a(t,s) 
 [e^{a(t,\tau)} -\kappa I]
 ( f(\tau) - f(t)) d\tau  
 -  U(t,t+\epsilon) f(t) + U(t,s)f(t), 
\end{array}  \end{equation}
where $\partial_{\tau} U(t,\tau) =  -A(\tau)  U(t,\tau)$ is utilized.
The last identity is obtained by applying $A(t) =  (I+\kappa U(s,t))\partial_{t} a(t,s)$. 
The H\"older continuity and Eq.~\eqref{leiq} lead to the strong convergence of the right hand of Eq.~\eqref{conveq}:
\[ \begin{array}{ll} 
 A(t)  \int_s^{t+\epsilon} [e^{a(t,\tau)} - \kappa I] f(\tau) d\tau   \vspace{1.5mm}\\
  \to
 \quad  \int_s^t (I +\kappa U(s,t)) (\partial_{t} a(t,s)) [e^{a(t,\tau)} -\kappa I]  (f(\tau) - f(t)) d\tau 
  + ( U(t,s) - I) f(t)
 \end{array}  \]
(due to $\epsilon \to 0$) for $f \in L^1(0,T;X)$. 
$A(t)$ is assumed to be an infinitesimal generator, so that $A(t)$ is a closed operator from $Y$ to $X$.
It follows that
\[   \begin{array}{ll}
  \int_s^{t} [e^{a(t,\tau)} -\kappa I] f(\tau) d\tau  \in Y
\end{array}  \]
and
\[ \begin{array}{ll} 
 A(t)  \int_s^{t} [e^{a(t,\tau)} -\kappa I] f(\tau) d\tau   \vspace{1.5mm}\\
  =
  \int_s^t (I + \kappa U(s,t)) (\partial_{t} a(t,s)) [e^{a(t,\tau)} -\kappa I]  (f(\tau) - f(t)) d\tau 
  + ( U(t,s) - I) f(t)  \in X.
 \end{array}  \]
The right hand side of this equation is strongly continuous on $[-T,T]$.
Consequently
  \[ \begin{array}{ll} 
  \partial_{t}  \int_s^{t+\epsilon} [e^{a(t,\tau)} -\kappa I] f(\tau) d\tau \vspace{1.5mm}\\
  =  [ e^{a(t,t+\epsilon)} -\kappa I ] f(t+\epsilon) +   \int_s^{t+\epsilon}  (\partial_{t} a(t,\tau)) e^{a(t,\tau)}  f(\tau) d\tau \vspace{1.5mm}\\
  \to  \quad f(t) + \int_s^{t} (I+\kappa U(\tau,t))^{-1}  A(t)  (U(t,\tau)+\kappa I)  f(\tau) d\tau \vspace{1.5mm}\\
 \qquad  =   f(t) +  \int_s^{t} A(t) U(t,\tau)   f(\tau) d\tau \vspace{1.5mm}\\
 \qquad  =   f(t) + A(t)  \int_s^{t} [e^{a(t,\tau)} -\kappa I] f(\tau) d\tau.
 \end{array}  \]
It is seen that $ \int_s^{t} [e^{a(t,\tau)} -\kappa I] f(\tau) d\tau$ satisfies Eq.~\eqref{origih},  and that it is sufficient to assume $f \in L^1(0,T;X)$ as H\"older continuous. \quad $\square$ \\ 

More simply, the unique solvability of non-autonomous case can be regarded in the context of decomposing the mild solution (for this terminology, see \cite{83pazy}).

\begin{corollary} \label{cor-inh}
Let $f \in L^1(-T,T;X)$ be locally H\"older continuous on $[-T,T]$.
For $u_s \in X$ there exists a unique solution $u(\cdot) \in C([-T,T];X)$ for \eqref{origih} such that
\[ \begin{array}{ll} 
 u(t)  =   \left[  \left( {\displaystyle \sum_{n=0}^{\infty} }  \frac{ a(t,s)^n}{n!}   \right) u_s 
 + \int_s^t    {\displaystyle \sum_{n=0}^{\infty}} \frac{ a(t,s)^n}{n!}     f(\tau) d\tau \right] 
  + \kappa I  \left[ u_s - \int_s^t   f(\tau) d\tau \right], \\
\end{array} \]
using a certain complex number $\kappa$.
\end{corollary}

{\bf Proof ~}
The representation is regarded as
\begin{equation}  \begin{array}{ll}
 u(t)  =   U(t,s) u_s + \int_s^tU(t,\tau) f(\tau) d\tau \vspace{1.5mm}\\
 =   ( e^{a(t,s)} - \kappa I )  u_s + \int_s^t  (  e^{a(t,\tau)} - \kappa I )  f(\tau) d\tau  \vspace{1.5mm}\\
 =   \left[ e^{a(t,s)} u_s + \int_s^t    e^{a(t,\tau)}   f(\tau) d\tau \right]  - \kappa I  \left[ u_s - \int_s^t   f(\tau) d\tau \right],
\end{array}  \end{equation}
The former part in the parenthesis is the mild solution\index{mild solution} of $\partial_t u(t)  = a(t,s)  u(t) + f(t)$, and the latter part in another parenthesis is the mild solution of $\partial_t u(t)  = f(t) $.
The unique existence of mild solution for the former part is valid for H\"older continuous $f \in L^1(0,T;X)$, and that for the latter part is valid for any $f \in L^1(-T,T;X)$. \quad $\square$ \\ 

Corollary \ref{cor-inh} shows the meaning of introducing the alternative infinitesimal generator.
This result should be compared to the standard theory of evolution equations in which the inhomogeneous term $f$ is assumed to be continuous on $[-T,T]$.
Consequently, in a purely abstract framework, the maximal regularity effect \cite{01pruess,16arendt} is found in the solutions of renormalized evolution equations. 
In this sense, the alternative infinitesimal generator brings about the analytic semigroup theory for non-parabolic evolution equations.

 \section{Relativistic formulation of abstract evolution equations}
\subsection{Formalism}
The relativistic formulation\index{relativistic formulation} of abstract evolution equations \cite{18iwata-2} is introduced to establish an abstract version of the Cole-Hopf transform in Banach spaces and to explain the nonlinear relation between the evolution operator and its infinitesimal generator \cite{18iwata-1}. 
The relativistic formulation is introduced for changing the evolution direction, which is necessary to justify the generalized Cole-Hopf transform.

In this paper the logarithmic representation of infinitesimal generator is utilized to formulate the relativistic form of abstract evolution equations.
Here the terminology ``relativistic" is used in the sense that there is no especially dominant direction.
In particular the role of $t$-direction (time direction) is not the absolute direction being compared to the other directions: $x$, $y$, and $z$-directions (spatial directions) in the standard notation.
While the relativistic treatment is associated with the equally-valid time-reversal and spatial-reversal symmetries, here the relativistic form to the generalized framework \eqref{replce} is introduced without assuming the invertible property of evolution operators.

Let the standard space-time variables $(t,x,y,z)$ be denoted by $(x^0,x^1,x^2,x^3)$ receptively.
It is further possible to generalize space-time variables to $(x^0,x^1,x^2,x^3, \cdots, x^{n})$ being valid to general $(n+1)$-dimensional space-time. 
In spite of the standard treatment of abstract evolution equations, the direction of evolution does not necessarily mean time-variable $t=x^0$ in the relativistic formulation of the abstract evolution equations.
Consequently, the equal treatment of any direction and the introduction of multi-dimension is naturally realized by the relativistic formulation.

\begin{definition}[Relativistic form]
For an evolution family of operators $\{U(x^i,\xi^i)\}_{-L \le x^i,\xi^i \le L}$ in a Banach space $X_i$, let  $K(x^i):Y_i \to X_i$ be the pre-infinitesimal generator of $U(x^i,\xi^i)$, where $Y_i$ is a dense subspace of $X_i$.
The relativistic form of abstract evolution equations is defined as
\begin{equation} \label{k-eq1}  \begin{array}{ll} 
\partial_{x^i} U(x^i, \xi^i)~u({\xi^i})  = K(x^i) U(x^i,\xi^i) ~u({\xi^i}) , \vspace{2.5mm} \\
 u({\xi^i}) = u_{\xi^i}
 \end{array}  \end{equation}
in $X_i$, where $X_i$ is a functional space consisting of functions with variables $x^j$ with $0 \le j \le n$ skipping only $j=i$.
Consequently, the unknown function is represented by $u(x^i)  =  U(x^i,\xi^i) ~u_{\xi^i}$ for a given initial value $u_{\xi^i} \in X_i$.
\end{definition}

Let $\partial_{x^i} U(x^i, \xi^i)~u_{\xi^i} = K(x^i) U(x^i,\xi^i) ~u_{\xi^i}$ evolving for $i$-direction be represented by $\partial_{x^k} V(x^k, \xi^k)~v_{\xi^k} =  {\mathcal K}(x^k) V(x^k,\xi^k) ~v_{\xi^k}$ in a certain direction $k$.
For $k \ne i$, let us begin with the abstract Cauchy problem
\begin{equation} \label{k-eq2}  \begin{array}{ll} 
\partial_{x^k} V(x^k, \xi^k)~v({\xi^k}) = {\mathcal K}(x^k) V(x^k,\xi^k) ~v({\xi^k}), \vspace{2.5mm} \\
 v({\xi^k}) = v_{\xi^k}
 \end{array}  \end{equation}
in $X_k$.
It is remarkable that even if the evolution operator $U(x^i, \xi^i)$ and its infinitesimal generator exist, $V(x^k, \xi^k)$ and its infinitesimal generator do not necessarily exist.
Those existence should be individually examined for each direction.
If $U(x^i, \xi^i)$, $V(x^k, \xi^k)$ and those infinitesimal generators exist, $u(x^i)$ in Eq.~(\ref{k-eq1}) and $v(x^k)$ in Eq.~(\ref{k-eq2}) satisfy the same evolution equation, where the detailed conditions such as initial and boundary conditions can be different depending on the settings of $X_i$ and $X_k$.  
For the purpose of introducing the relativistic form with a significance, it is necessary to clarify
\begin{itemize}
\item the well-defined (pre-)infinitesimal generator of $V(x^k, \xi^k)$
\item the existence of $V(x^k, \xi^k)$ (or the corresponding regularized evolution operator)
\end{itemize}
to an unknown direction $k$.
The second one automatically follows if the first one is established.
Otherwise Eq.\ref{k-eq2} cannot be regarded as the abstract evolution equations.
This issue is examined in generalizing the Cole-Hopf transform.

The propagation of singularity should be different if the evolution direction is different. 
For Eqs.~(\ref{k-eq1}) and (\ref{k-eq2}), the evolution direction is not limited to $x^0$.
This gives a reason why the formulation shown in Eq.~(\ref{k-eq1}) is called the relativistic form of abstract evolution equations. 
It means that if invertible evolution operator is obtained for one direction, the evolution operator for the other direction is not necessarily be the invertible.
Here is a reason why it is useful to introduce a relativistic form based on the generalized logarithmic representation (Cf.~Corollary~\ref{cor1}).

One utility of considering the evolution towards spatial direction is to explain and generalize the Cole-Hopf transform.
For this purpose, it is necessary to realize the logarithmic representation of the infinitesimal generators defined in the relativistic form of the abstract evolution equations.  
That is, for a significant introduction of the relativistic form, it should be introduced together with the logarithmic representation.
The condition to obtain the logarithmic representation is stated as follows.

\begin{theorem}[Relativistic form of logarithmic representation] \label{thm-rel}
Let $i$ denote any direction satisfying $0 \le i \le n$.
Let $x^i$ and $\xi^i$ satisfy $-L \le x^i,\xi^i \le L$, and $Y_i$ be a dense subspace of a Banach space $X_i$.
A two-parameter evolution family of operators $\{U(x^i,\xi^i)\}_{-L \le x^i,\xi^i \le L}$ satisfying Eq.~(\ref{sg1}) is assumed to exist in a Banach space $X_i$ (i.e., the inverse of $U(x^i,\xi^i)$ is not assumed).
Under the existence of the pre-infinitesimal generator $K(x^i):Y_i \to X_i$ of $U(x^i,\xi^i)$ for the $x^i$ direction, let $U(x^i,\xi^i)$ and $K(x^i)$ commute.
The logarithmic representation of infinitesimal generator is obtained; there exists a certain complex number $\kappa \ne 0$ such that
\begin{equation} \label{logex2} \begin{array}{ll}
    K(x^i) ~ u_{\xi^i}
    = (  I - \kappa e^{-a(x^i,\xi^i)} )^{-1} ~ \partial_{x^i} {\rm Log} ~ (U(x^i,\xi^i) + \kappa I) ~ u_{\xi^i}
 =  (  I - \kappa e^{-a(x^i,\xi^i)} )^{-1} ~ \partial_{x^i} a(x^i,\xi^i) ~ u_{\xi^i},
\end{array} \end{equation}
where $u_{\xi}$ is an element in $Y_i$, $\kappa$ is taken from the resolvent set of $U(x^i,\xi^i)$, and $a(x^i,\xi^i) = {\rm Log} ~ (U(x^i,\xi^i) + \kappa I)$.
Note that $U(x^i,\xi^i)$ is not assumed to be invertible.
\end{theorem}
{\bf Proof ~} 
Different from the proof of Theorem~\ref{thm1}, here the similar statement is proved without assuming the invertible property of $U(x^i,\xi^i)$.
The key point is that $  (  I - \kappa e^{-a(x^i,\xi^i)} )^{-1}$ exists for a certain $\kappa \in {\mathbb C}$, even if $ U(x^i,\xi^i)^{-1} $ does not exist.
In particular, the obtained representation is more generally compared to the one obtained in Ref.~\cite{17iwata-1}.
For any $U(x^i,\xi^i)$, operators $ {\rm Log} ~ (U(x^i,\xi^i) + \kappa I)$ and $ {\rm Log} ~ (U(x^i+h,\xi^i) + \kappa I)$ are well defined for a certain $\kappa$.
The $x^i$-differential in a weak sense is formally written by
\begin{equation} \label{difference0-2} \begin{array} {ll} 
\mathop{\rm wlim}\limits_{h \to 0}  \frac{1}{h} \{ {\rm Log} ~(U(x^i+h,\xi^i)+\kappa I) - {\rm Log} ~(U(x^i,\xi^i)+ \kappa I) \}   \vspace{1.5mm} \\
 = \mathop{\rm wlim}\limits_{h \to 0}  \frac{1}{2 \pi i}
 \int_{\Gamma} {\rm Log} \lambda   
 ~ \{ (\lambda I - U(x^i+h,\xi^i)-\kappa I )^{-1} \frac{U(x^i+h,\xi^i)-U(x^i,\xi^i)}{h} (\lambda I - U(x^i,\xi^i) - \kappa I )^{-1} \} d \lambda   
\end{array} \end{equation} 
where $\Gamma$, which is taken independent of $x^i$, $\xi^i$ and $h$ for a sufficiently large certain $\kappa$, denotes a circle in the resolvent set of both $U(t,s)+ \kappa I$ and $U(t+h,s)+\kappa I$.

The discussion, which is the same as that shown in Theorem~\ref{thm1}, leads to
\[ \begin{array} {ll} 
\partial_t {\rm Log} (U(x^i,\xi^i) + \kappa I) ~ u \vspace{1.5mm}\\
=  \frac{1}{2 \pi i} \int_{\Gamma} d\lambda    
 ~ \left[ ( {\rm Log} \lambda )  (\lambda I - U(x^i,\xi^i) - \kappa I )^{-1} 
 ~\mathop{\rm wlim}\limits_{h \to 0} \left ( \frac{U(x^i+h,\xi^i)-U(x^i,\xi^i)}{h} \right)
~ (\lambda I  - U(x^i,\xi^i)  - \kappa I )^{-1}    \right] ~u \\
\end{array} \] 
for $u \in Y_i$.
Because it is allowed to interchange $K(x^i)$ with $U(x^i,\xi^i)$,
\[  \begin{array}{ll}
\partial_t {\rm Log} (U(x^i,\xi^i) + \kappa I) ~ u 
 = \frac{1}{2 \pi i} \int_{\Gamma} ({\rm Log} \lambda) ~(\lambda I -U(x^i,\xi^i) - \kappa I)^{-2} ~ U(t,s) ~ d \lambda ~ K(x^i) ~ u
\end{array} \]
for $u \in Y_i$, where $\mathop{\rm wlim}\limits_{h \to 0} \left ({U(x^i+h,\xi^i)-U(x^i,\xi^i)} \right)/h$ means the pre-infinitesimal generator $K(x^i)$ itself.
A part of the right hand side is calculated as
\[ \begin{array}{ll}
 \quad  \frac{1}{2 \pi i} \int_{\Gamma}~ ({\rm Log} \lambda) ~(\lambda I -U(x^i,\xi^i) - \kappa I)^{-2} U(x^i,\xi^i) ~ d \lambda  
  = I- \kappa (U(t,s)+ \kappa I)^{-1},  
\end{array} \] 
due to the integration by parts, where the details of procedure is essentially the same as Ref.~\cite{17iwata-1}.
It leads to
\[ \begin{array}{ll}
 K(x^i) ~ u_{\xi^i} =
  (  I - \kappa e^{-a(x^i,\xi^i)} )^{-1} ~  \partial_{x^i} {\rm Log} ~ (U(x^i,\xi^i) + \kappa I) ~ u_{\xi^i},
\end{array} \]
for $u_{\xi^i} \in Y_i$.
It is notable that $ (U(x^i,\xi^i) + \kappa I)^{-1}$ is always well defined for any $\kappa$ taken from the resolvent set of $U(x^i,\xi^i)$, even if $ U(\xi^i,x^i) = U(x^i,\xi^i)^{-1}$ does not exist.
\quad $\square$ \\

Under the existence of logarithmic representation for $K(x^i)$, the related concepts such as
\begin{enumerate}
\item alternative infinitesimal generator: $\partial_{x^i} a(x^i,\xi^i) =  \partial_{x^i} {\rm Log} ( U(x^i,\xi^i) + \kappa I)$
\item regularized evolution operator: $e^{a(x^i,\xi^i)} = U(x^i,\xi^i) + \kappa I$
\item renormalized abstract evolution equation: $\partial_{x^i} {\tilde u} = [\partial_{x^i} a(x^i,\xi^i) ] {\tilde u}$
\end{enumerate}
are similarly well-defined in the relativistic framework.

\subsection{Generalization of  the Cole-Hopf transform} 
Now it is ready for establishing the general version of the Cole-Hopf transform.
It corresponds to an application example of relativistic formulation is provided.
The Cole-Hopf transform \cite{15bateman, 48burgers, 51cole, 06forsyth, 50hopf} is a concept bridging the linearity and the nonlinearity.
In the following, such a linear-nonlinear conversion relation is found within the relation between the infinitesimal generators and the generated semigroups.

For $t \in {\mathbb R}_+$ and $x \in {\mathbb R}$ the Cole-Hopf transform\index{Cole-Hopf transform} reads
\begin{equation}  \begin{array}{ll}
 \psi(t,x) = -2 \mu^{-1/2} ~ \partial_x \log u(t,x), 
\end{array} \end{equation}
where $u(t,x)$ denotes the solution of linear equation, and $\psi(t,x)$ is the solution of transformed nonlinear equation.
On the other hand, e.g., for $t \in {\mathbb R}$ and $x \in {\mathbb R}^n$, the logarithmic representation of infinitesimal generator 
\begin{equation} \label{logex2} \begin{array}{ll}
 K(x^i) ~ u_{\xi^i} =
  (I+ \kappa U(\xi^i,x^i))~ \partial_{x^i} {\rm Log} ~ (U(x^i,\xi^i) + \kappa I) ~ u_{\xi^i},
\end{array} \end{equation}
has been obtained in the abstract framework, where $U(x^i,\xi^i)$ denotes the evolution operator, and $ K(x^i) $ is its infinitesimal generator.  
By taking a specific case with $\kappa = 0$, the similarity between them is clear.
That is, the process of obtaining infinitesimal generators from evolution operators is expected to be related to the emergence of nonlinearity.

Based on the logarithmic representation of infinitesimal generators obtained in Banach spaces, the Cole-Hopf transform is generalized in the following sense:
\begin{itemize}
\item the linear equation is not necessarily the heat equation;
\item the spatial dimension of the equations is not limited to 1;
\item the variable in the transform is not limited to a spatial variable $x$;
\end{itemize}
where, in order to realize these features, the relativistic formulation of abstract evolution equation is newly introduced.
Since the logarithmic representation shows a relation between an evolution operator and its infinitesimal generator, the correspondence to the Cole-Hopf transform means a possible appearance of nonlinearity in the process of defining an infinitesimal generator from the evolution operator.

The next theorem follows.

\begin{theorem}[Generalization of the Cole-Hopf transform]    \label{thm3}
Let $i$ be an integer satisfying $0 \le i \le n$, and $Y$ be a dense subspace of Banach space $X$.
Let an invertible evolution family $\{ U(x^i, \xi^i)\}_{0 \le x^i,\xi^i \le L}$ be generated by $A(x^i)$ for $0 \le x^i,\xi^i \le L$ in a Banach space $X$.
$U(x^i, \xi^i)$ and $A(x^i)$ are assumed to commute.
For any $u_{\xi^i} \in Y \subset X$, the logarithmic representation
\begin{equation} \label{spatrep3} \begin{array}{ll}
 A(x^i)  U(x^i, \xi^i) u_{\xi^i}  
  =  e^{a(x^i,\xi^i)}~ [ \partial_{x^i} {\rm Log} ~ ( U(x^i, \xi^i) + \kappa I) ]  u_{\xi^i} ,
\end{array} \end{equation}
is the generalization of the Cole-Hopf transform, where the logarithmic representation is obtained in a general Banach space framework, $\kappa \ne 0$ is a complex number and , where $a(x^i, \xi^i) = {\rm Log} (U(x^i, \xi^i) +  \kappa I)$.
In particular, if $(U(x^i, \xi^i) u_{\xi^i})^{-1}$ exists for a given interval $0 \le x^i,\xi^i \le L$, its normalization
\begin{equation} \label{spatrep5} \begin{array}{ll}
 A(x^i)  
  =  ( I - \kappa e^{-a(x^i, \xi^i)})~ [ \partial_{x^i} {\rm Log} ~ ( U(x^i, \xi^i) + \kappa I) ]
\end{array} \end{equation}
defined in $X$ corresponds to $\psi(t,x)$.
\end{theorem}

{\bf Proof ~} 
The proof consists of five steps. \\

{\bf [1st step: formulation] ~}
It is necessary to recognize the evolution direction of the heat equation\index{heat equation} as $x$, because the derivative on the spatial direction $x$ is considered in the Cole-Hopf transform.
The Cole-Hopf transform acts on one-dimensional heat equation
\begin{equation} \label{heat} \begin{array}{ll}
\partial_x^2 u(t,x) -\mu^{1/2}  \partial_t u(t,x) = 0, \quad t \in (0, \infty) , ~ x \in (-L, L), \vspace{2.5mm}\\
u(t,-L) = u(t,L)  = 0,  \quad t \in (0, \infty), \vspace{2.5mm}\\
u(0,x) = u_0(x),  \quad x \in (-L, L),
\end{array} \end{equation}
where $\mu$ is a real positive number, and the hypoelliptic\index{hypoelliptic} property of parabolic evolution equation is true.
The first equation of (\ref{heat}) is hypoelliptic; for an open set ${\mathcal U} \subset (-\infty,\infty) \times  (-L,L)$, $u \in C^{\infty}({\mathcal U})$ follows from $(  \partial_x^2 - \mu^{1/2} \partial_t) u \in C^{\infty}({\mathcal U})$.
Equation~(\ref{heat}) is well-posed in $C^{\infty}(0,\infty) \times C^{\infty}(-L,L)$, so that $\mu^{1/2}  \partial_t$ is the infinitesimal generator in $C^{\infty}(-L,L)$.
The spaces $C^{1}(-L,L)$ and  $C^{\infty}(-L,L)$ are dense in $L^{p}(-L,L)$ with $1 \le p < \infty$.
The solution is represented by
\[ \begin{array}{ll}
 u(t,x) = U(t) u_0, \ 
\end{array} \]
where $U(t)$ is a semigroup generated by $\mu^{-1/2}  \partial_x^2$ under the Dirichlet-zero boundary condition.

By changing the evolution direction from $t$ to $x$, the heat equation
\begin{equation} \label{heat2} \begin{array}{ll}
\partial_x^2 u(t,x) -\mu^{1/2}  \partial_t u(t,x) = 0, \quad t \in (0, \infty) , ~ x \in (-L, L), \vspace{2.5mm}\\
u(t,-L) = v_0(t), \quad \partial_x u(t,-L) = v_1(t),  \quad t \in (0, \infty), \vspace{2.5mm}\\
u(0,x) = 0,  \quad x \in (-L, L),
\end{array} \end{equation}
is considered for $x$-direction, where $v_0(t)$ and $v_1(t)$ are given initial functions.
To establish the existence of semigroup for the $x$-direction, it is sufficient to consider the generation of semigroup in $L^2(-\infty, \infty)$ by generalizing $t$-interval from $(0, \infty)$ to $(-\infty, \infty)$.
The Fourier transform leads to
\begin{equation}  \begin{array}{ll}
\partial_x^2 \tilde{u} - i \mu^{1/2} \omega \tilde{u}  = 0 \vspace{2.5mm}\\
\tilde{u} (\omega,0)  =  {\tilde v}_0(\omega)  , \quad  
\partial_x \tilde{u} (\omega,0) =  {\tilde v}_1(\omega) , 
\end{array} \label{eq-add} \end{equation}
where $\omega$ is a real number.
Indeed, the following transforms
\[  \begin{array}{ll}
u (t,x) = \frac{1}{2 \pi} \int_{-\infty}^{\infty}  \tilde{u}(\omega,x)  e^{ i \omega t} d \omega, \vspace{1.5mm} \\
u (t, -L) = \frac{1}{2 \pi}  \int_{-\infty}^{\infty} {\tilde u} (\omega, -L) e^{ i \omega t} d \omega    \vspace{1.5mm}, \\
\partial_x  u (t, -L) = \frac{1}{2 \pi}  \int_{-\infty}^{\infty}   \partial_x {\tilde u} (\omega, -L)  e^{i \omega t} d \omega  
\end{array} \]
are implemented.
By solving the characteristic equation $\lambda^2 -i \mu^{1/2} \omega = 0$, the Fourier transformed solution of (\ref{heat2}) is
\[  \begin{array}{ll}
\tilde{u}(\omega, x) 
 = \frac{ {\tilde v}_0(\omega)- i(- i \mu^{1/2}\omega )^{-1/2} {\tilde v}_1(\omega) }{2}  e^{ + ( i \mu^{1/2} \omega)^{1/2}  x}  
+  \frac{ {\tilde v}_0(\omega)+i  (- i \mu^{1/2}\omega)^{-1/2} {\tilde v}_1(\omega) }{2} e^{ - ( i \mu^{1/2} \omega)^{1/2}  x },
\end{array} \]
where
\[  \begin{array}{ll}
\tilde{u}(\omega, 0) 
 = \frac{ {\tilde v}_0(\omega)- i(- i \mu^{1/2}\omega)^{1/2} {\tilde v}_1(\omega) }{2}  +  \frac{ {\tilde v}_0(\omega)+i  (- i \mu^{1/2}\omega)^{1/2}  {\tilde v}_1(\omega)}{2} = {\tilde v}_0(\omega),
  \vspace{2.5mm} \\
\partial_x \tilde{u}(\omega, 0) 
 =  ( i \mu^{1/2} \omega)^{1/2} \left( \frac{ {\tilde v}_0(\omega)- i(- i \mu^{1/2}\omega)^{-1/2} {\tilde v}_1(\omega) }{2}  -  \frac{ {\tilde v}_0(\omega)+i  (- i \mu^{1/2}\omega)^{-1/2} {\tilde v}_1(\omega) }{2}\right) =  {\tilde v}_1(\omega).
\end{array} \]
Meanwhile, based on the relativistic treatment, one-dimensional heat equation
\[ \begin{array}{ll}
\partial_x^2 u -\mu^{1/2}  \partial_t u = 0
\end{array} \]
is written as
\[ \begin{array}{ll}
\partial_x 
\left( \begin{array}{ll}
 u \\
 v
 \end{array}  \right) 
  -
\left( \begin{array}{cc}
0 & I \\f
 \mu^{1/2}  \partial_t   & 0
\end{array} \right)
\left( \begin{array}{ll}
 u \\
 v
 \end{array}  \right)
 = 0. \vspace{1.5mm}
\end{array} \]
Let a linear operator ${\mathcal A}$ be defined by
\[ \begin{array}{ll}
{ \mathcal A} = 
\left( \begin{array}{cc}
0 & I \\
\mu^{1/2}  \partial_t   & 0
\end{array} \right)
 \end{array} \]
in $L^2(-\infty,\infty) \times L^2(-\infty,\infty)$, and the domain space of ${\mathcal A}$ be
\[ \begin{array}{ll} 
D({\mathcal A})= H^1(-\infty,\infty) \times L^2(-\infty,\infty), 
\end{array} \]
where $H^1(-\infty,\infty) = \{ u(t) \in L^2(-\infty,\infty);   ~u'(t) \in L^2(-\infty,\infty), ~ u(0) = 0  \}$ is a Sobolev space\index{Sobolev space}.
The Fourier transform means that the diagonalization of ${\mathcal A} $ is equal to
\[ \begin{array}{ll}
\tilde{\mathcal A} = 
\left( \begin{array}{cc}
(\mu^{1/2}  \partial_t )^{1/2}  & 0 \\
0  & - ( \mu^{1/2}  \partial_t )^{-1/2} 
\end{array} \right).
 \end{array} \]
In this context the master equation of the problem (\ref{heat2}) is reduced to the abstract evolution equation
\begin{equation} \begin{array}{ll}
\partial_x 
\left( \begin{array}{ll}
 {\tilde u} \\
 {\tilde v}
 \end{array}  \right) 
  -
\tilde{\mathcal A}
\left( \begin{array}{ll}
 {\tilde u} \\
 {\tilde v}
 \end{array}  \right)
 = 0
\end{array} \label{abs} \end{equation}
in $L^2(-\infty,\infty) \times L^2(-\infty,\infty)$, where 
\[  \begin{array}{ll}
\tilde{u}=
 \frac{u - i(- i \mu^{1/2}\omega )^{-1/2} v }{2}  
\end{array} \]
and
\[  \begin{array}{ll}
\tilde{v}=
 \frac{ u +i  (- i \mu^{1/2}\omega)^{-1/2} v }{2}.
\end{array} \]
It suggests that the evolution operator of Eq.~(\ref{heat2}) is generated by
\begin{equation}  \begin{array}{ll}
\pm   (\mu^{1/2} \partial_t)^{1/2}, 
\end{array} \end{equation}
so that it is sufficient to show $\pm   (  \mu^{1/2} \partial_t)^{1/2}$ as the infinitesimal generator.
Note that the operator $\tilde{\mathcal A}$ is not necessarily a generator of analytic semigroup, because the propagation of singularity should be different if the evolution direction is different.
Consequently the existence of semigroup for (\ref{abs}) in the $x$-direction is reduced to show $\pm (\mu^{1/2} \partial_t)^{1/2}$ as the infinitesimal generator in $L^2(-\infty, \infty)$.
In the following, the property of $\mu^{1/2} \partial_t$ is discussed in the second step, and the fractional power\index{fractional power} of $ (  \mu^{1/2} \partial_t)^{1/2}$ is studied in the third step.  \\

{\bf [2nd step: first order differential operator] ~}
The following lemma is proved in this step.

\begin{lemma} \label{lem-m1}
The operator $\mu^{1/2} \partial_t$ with the domain $H^1(-\infty, \infty)$ is the infinitesimal generator in $L^2(-\infty, \infty)$.
\label{lemma1} \end{lemma}

{\bf Proof of Lemma \ref{lem-m1}.} 
Let $\lambda$ be a complex number satisfying ${\rm Re }\lambda >0$.
First the existence of $(\lambda - \mu^{1/2} \partial_t)^{-1}$ is examined.
Let $f(t)$ be included in $L^2 (-\infty,\infty)$. 
Because
\begin{equation} \begin{array}{ll}
(\lambda - \mu^{1/2} \partial_t) u = f
\end{array} \label{iheq} \end{equation}
in one-dimensional interval $(-\infty,\infty)$ is a first-order ordinary differential equation with a constant coefficient, and the global-in-$t$ solution necessarily exists for a given $u(0)=u_0 \in {\mathbb C}$.
That is, $\lambda / \mu^{1/2}$ is included in the resolvent set of $\partial_t$ for an arbitrary complex number $\lambda$, so that $(\lambda - \mu^{1/2} \partial_t)^{-1}$ is concluded to be well-defined in $H^1(-\infty,\infty)$.

Second the resolvent operator $(\lambda - \mu^{1/2} \partial_t)^{-1}$ is estimated from the above.
Since $\lambda / \mu^{1/2}$ is included in the resolvent set of $\partial_t$,
it is readily seen that $(\lambda - \mu^{1/2} \partial_t)^{-1}$ is a bounded operator on $L^2(-\infty, \infty)$.
More precisely, let us consider Eq.~(\ref{iheq}) being equivalent to
\begin{equation}  \begin{array}{ll}
 \partial_t u(t) =   (\lambda / \mu^{1/2}) u(t) - f(t) / \mu^{1/2}.
\end{array}  \label{sim1} \end{equation}
If the inhomogeneous term satisfies $f(t) \in L^2 (-\infty,\infty)$, 
\begin{equation} \begin{array}{ll}
u(t) = 
- \frac{1}{\mu^{1/2}} \int_t^{\infty}  \exp \left( \frac{- \lambda }{ \mu^{1/2} }  (s-t) \right) f(s)  ds
\end{array} \label{intrep} \end{equation}
satisfies Eq.~(\ref{sim1}).
According to the Schwarz inequality,
\[ \begin{array}{ll}
\int_{-\infty}^{\infty} |u(t)|^2 dt  
 = \int_{-\infty}^{\infty} \left| \frac{1}{\mu^{1/2}} \int_t^{\infty}  \exp \left( \frac{- \lambda }{ \mu^{1/2} }  (s-t) \right) f(s)  ds  \right|^2 dt      \vspace{2.5mm}  \\
\quad \le \frac{1}{\mu} \int_{-\infty}^{\infty} \left\{ \int_t^{\infty} 
  \exp \left( \frac{- {\rm Re} \lambda }{ 2 \mu^{1/2} }  (s-t) \right) \right. 
  \left. \exp \left( \frac{- {\rm Re} \lambda }{ 2 \mu^{1/2} }  (s-t) \right) | f(s) |  ds  \right\}^2 dt      \vspace{2.5mm}  \\
\quad \le \frac{1}{\mu} \int_{-\infty}^{\infty}   \int_t^{\infty} 
  \exp \left( \frac{-  {\rm Re} \lambda }{ \mu^{1/2} }  (s-t) \right) ds 
 ~ \int_t^{\infty}  \exp \left( \frac{-  {\rm Re} \lambda }{ \mu^{1/2} }  (s-t) \right) | f(s) |^2  ds ~ dt
\end{array} \]
is obtained, and the equality
\[ \begin{array}{ll}
\int_{t}^{\infty}  \exp \left( \frac{-  {\rm Re} \lambda }{ \mu^{1/2} }  (s-t) \right) ds 
= \int_{0}^{\infty}  \exp \left( \frac{- {\rm Re} \lambda }{ \mu^{1/2} } s \right) ds 
=  \frac{ \mu^{1/2} }{  {\rm Re} \lambda }
\end{array} \]
is positive if ${\rm Re}\lambda >0$ is satisfied.
Its application leads to
\[ \begin{array}{ll}
\int_{-\infty}^{\infty} |u(t)|^2 dt    
 \le \frac{1}{\mu} \frac{ \mu^{1/2} }{  {\rm Re} \lambda }   \int_{-\infty}^{\infty}  
  \int_t^{\infty}  \exp \left( \frac{-  {\rm Re} \lambda }{ \mu^{1/2} }  (s-t) \right) | f(s) |^2  ds ~ dt    \vspace{2.5mm}  \\
\quad = \frac{1}{\mu} \frac{ \mu^{1/2} }{  {\rm Re} \lambda }   \int_{-\infty}^{\infty}  
  \int_{-\infty}^{s}  \exp \left( \frac{-  {\rm Re} \lambda }{ \mu^{1/2} }  (s-t) \right) ~ dt  ~  | f(s) |^2  ~ ds.
\end{array} \]
Further application of the equality
\[ \begin{array}{ll}
\int_{-\infty}^{s}  \exp \left( \frac{-  {\rm Re} \lambda }{ \mu^{1/2} }  (s-t) \right) dt
= \int_{-\infty}^{0}  \exp \left( \frac{  {\rm Re} \lambda }{ \mu^{1/2} } t \right) dt 
=  \frac{ \mu^{1/2} }{  {\rm Re} \lambda }
\end{array} \]
results in
\[ \begin{array}{ll}
\int_{-\infty}^{\infty} |u(t)|^2 dt    \vspace{2.5mm}  
 \le  \frac{1}{\mu} \frac{ \mu^{1/2} }{  {\rm Re} \lambda }   \int_{-\infty}^{\infty}  
 \frac{ \mu^{1/2} }{ {\rm Re} \lambda } ~  | f(s) |^2  ~ ds  
 = \frac{ 1 }{  {\rm Re} \lambda^2 }   \int_{-\infty}^{\infty}   ~  | f(s) |^2  ~ ds,
\end{array} \]
for ${\rm Re}\lambda >0$, and therefore
\[ \begin{array}{ll} 
\| (\lambda I - \mu^{1/2} \partial_t)^{-1} f(t) \|_{L^2(-\infty,\infty)} 
\le \frac{1}{ {\rm Re} \lambda} \| f(t) \|_{L^2(-\infty,\infty)}
\end{array} \]
follows. 
That is, for ${\rm Re}\lambda >0$,
\begin{equation} \label{hyouka-1} \begin{array}{ll}
\|(\lambda I - \mu^{1/2} \partial_t)^{-1} \| \le 1/  {\rm Re} \lambda
\end{array} \end{equation}
is valid.
The surjective property of $\lambda I - \mu^{1/2} \partial_t$ is seen by the unique existence of solutions $u \in L^2(-\infty,\infty)$ for the Cauchy problem of Eq.~(\ref{sim1}).

A semigroup is generated by taking a subset of the complex plane as 
\[ \begin{array} {ll}
\Omega = \left\{ \lambda \in {\mathbb C}; ~ \lambda = {\bar \lambda}  \right\},
\end{array} \]
where $\Omega$ is included in the resolvent set of $\mu^{1/2} \partial_t$.
For $\lambda \in \Omega$, $(\lambda I - \mu^{1/2} \partial_t)^{-1}$ exists, and
\begin{equation} \begin{array}{ll}
\| (\lambda I - \mu^{1/2} \partial_t)^{-n} \| 
 \le 1/({\rm Re} \lambda)^{n}
\end{array} \label{resolves} \end{equation}
is obtained.
Consequently, according to the Lumer-Phillips theorem\index{Lumer-Phillips theorem}~\cite{61lumer} for the generation of quasi contraction semigroup, $\mu^{1/2} \partial_t$ is confirmed to be an infinitesimal generator in $L^2(-\infty,\infty)$, and the unique existence of global-in-$x$ weak solution follows. \quad {\bf [Q.E.D.: Lemma \ref{lem-m1}]}  \vspace{2.5mm}\\ 

The semigroup generated by $\mu^{1/2} \partial_t$ is represented by
\[ \begin{array}{ll}
\left(  \exp (h \mu^{1/2} \partial_t) u \right) (t) = u(t+\mu^{1/2}h), \quad -\infty < h < \infty,
\end{array} \]
so that the group is actually generated by $\mu^{1/2} \partial_t$.
Indeed, the similar estimate as Eq.~(\ref{resolves}) can be obtained for ${\rm Re} \lambda < 0$ with $(\lambda - \mu^{1/2})u=f$ in which the solution $u$ is represented by
\[ \begin{array}{ll}
u(t) = -\frac{1}{\mu^{1/2}} \int_{-\infty}^t \exp \left( \frac{-\lambda}{\mu^{1/2}} (s-t) \right) f(s) ds
\end{array} \]
that should be compared to Eq.~(\ref{intrep}). \\

{\bf [3rd step: fractional powers of operator] ~}
The following lemma is proved in this step.

\begin{lemma} \label{lem-m2}
For $0<\alpha <1$, the operator $(\mu^{1/2} \partial_t)^{\alpha}$ is the infinitesimal generator in $L^2(-\infty, \infty)$.
\label{thm1a}
\end{lemma}

 {\bf [Proof of Lemma \ref{lem-m2}].}  
According to Lemma~\ref{lemma1}, $\mu^{1/2} \partial_t$ is the infinitesimal generator in $L^2(-\infty,\infty)$.
For an infinitesimal generator $\mu^{1/2} \partial_t$ in $L^2(-\infty,\infty)$, let the one-parameter semigroup generated by $\mu^{1/2} \partial_t$ be denoted by ${\mathcal V}(x)$.
An infinitesimal generator $\mu^{1/2} \partial_t$ is a closed linear operator in $L^2(-\infty,\infty)$.
Its fractional power
\[ \begin{array} {ll}
(\mu^{1/2} \partial_t)^{\alpha}, \quad 0<\alpha<1
\end{array} \]
has been confirmed to be well-defined by S. Bochner~\cite{49bochner} and R.S. Phillips~\cite{52phillips} as the infinitesimal generator of semigroup (cf.~K. Yosida~\cite{{60yosida}}):
\[ \begin{array} {ll}
W(x) w_0 = \int_0^{\infty} {\mathcal V}(x) w_0 ~ d \gamma(\lambda),
\end{array} \]
for $w_0 \in L^2(-\infty, \infty)$, where $W(x)$ is the semigroup for $x$-direction.
The measure $d \gamma(\lambda) \ge 0$ is defined through the Laplace integral
\[ \begin{array} {ll}
\exp (-tk^{\alpha}) = \int_0^{\infty} \exp(- \lambda k) ~ d \gamma(\lambda),
\end{array} \]
where $t, k>0$ is satisfied. \quad  {\bf [Q.E.D.: Lemma \ref{lem-m2}]}    \vspace{2.5mm}\\

By taking $\alpha = 1/2$, $(\mu^{1/2} \partial_t)^{1/2}$ is confirmed to be an infinitesimal generator in $L^2(-\infty, \infty)$.
Because  $\Omega$ is included in the resolvent set of $-(\mu^{1/2} \partial_t)^{1/2}$, it is readily seen that $-(\mu^{1/2} \partial_t)^{1/2}$ is also an infinitesimal generator in $L^2(-\infty, \infty)$. \\

{\bf [4th step: abstract form of the Cole-Hopf transform] ~}
As in the original derivation of the Cole-Hopf transform, the solution of heat equation $w(t,x) = W(x,\xi) w_{\xi}(t)$ solved along the $x$-direction permits its logarithm function.
The abstract case of original Cole-Hopf transform is included in the description of the logarithmic representation (\ref{spatrep3}).
Indeed, let an invertible evolution family $\{ W(x, \xi)\}_{0 \le x,\xi \le L}$ be generated by ${\mathcal A}(x)$ for $0 \le x,\xi \le L$.
According to Lemma~\ref{thm1a} the logarithmic representation of relativistic form (\ref{logex2}) is obtained in this case as
\[  \begin{array}{ll}
 {\mathcal A}(x) w_{\xi} =  (I+ \kappa W(\xi, x))~ [ \partial_{x} {\rm Log} ~ (W(x,\xi) + \kappa I) ] w_{\xi}
\end{array} \]
and hence as
\[  \begin{array}{ll}
 {\mathcal A}(x)  W(x, \xi) w_{\xi} 
  =  (\kappa I + W(x, \xi))~ [ \partial_{x} {\rm Log} ~ (W(x,\xi) + \kappa I) ] w_{\xi},
\end{array} \]
using the commutation assumption.
The nonlinear Anzatz $-2 \mu^{-1/2} (\partial_x u(t,x)) ~ u(t,x)^{-1}$ of the Burgers' equation\index{Burgers' equation}
\begin{equation}  \label{ch-trans} \begin{array}{ll}
 \psi(t,x) = -2 \mu^{-1/2} ~ \partial_x \log u(t,x)
  = -2 \mu^{-1/2} (\partial_x u(t,x)) ~ u(t,x)^{-1}, 
\end{array} \end{equation}
is essentially represented by
\begin{equation} \label{spatrep4} \begin{array}{ll}
 -2 \mu^{-1/2}({\mathcal A}(x)  W(x, \xi)) ~  W(x, \xi)^{-1}   \vspace{1.5mm}\\
  = -2 \mu^{-1/2}  (\kappa I + W(x, \xi))~ [ \partial_{x} {\rm Log} ~ (W(x,\xi) + \kappa I) ] ~  W(x, \xi)^{-1}  \vspace{1.5mm}\\
  = -2 \mu^{-1/2}  (\kappa  W(\xi, x) + I)~ [ \partial_{x} {\rm Log} ~ (W(x,\xi) + \kappa I) ]
\end{array} \end{equation}
in the  abstract form.
The similarity between Eq.~(\ref{spatrep5}) and the standard definition of operator norm is clear.
In particular the evolution direction is generalized from $x$ to $x^i$ in Eq.~(\ref{spatrep5}). \\

{\bf [5th step: generalization property] ~}
Equation ~(\ref{spatrep3}) is the generalization of the Cole-Hopf transform.
According to the introduction of nonzero $\kappa$ in the abstract form, the applicability is significantly increased, so that the linear equation is not necessarily the heat equation.
According to the abstract nature of the logarithmic representation, linear and nonlinear equations are not necessarily considered in the one-spatial dimension.
According to the relativistic treatment, the transformed variable is not limited to the spatial variables. 
 \quad $\square$ \\

The generalized Cole-Hopf transform (\ref{spatrep3})  shows that the nonlinearity of semigroup can appear simply by altering the evolution direction under a suitable identification between the infinitesimal generator and the evolution operator.
In this sense, Eq.~(\ref{k-eq1}) is regarded as a local-in-$x^i$ linearized equation, if $U(x^i,\xi^i)$ is a nonlinear semigroup (semigroup related to the nonlinear equations).
Furthermore, the generalized Cole-Hopf transform (\ref{spatrep3})  suggests that
the relation between evolution operator and its infinitesimal generator corresponds essentially to the transform between linearity and nonlinearity. 
In the same context of generalizing Miura transform between the Korteweg-de-Vries and the modified Korteweg-de-Vries equations, the logarithmic representation is utilized~\cite{20iwata}. 
\vspace{3mm} \\.

\section{Algebraic structure of infinitesimal generators}
\subsection{$B(X)$-module}
The algebraic structure is studied based on the relativistic form of abstract equations.
The operator $a(x^i,\xi^i) = {\rm Log} ~ (U(x^i,\xi^i) + \kappa I)$ is bounded on $X$.
It follows that $e^{a(x^i,\xi^i)}$ is well-defined by the convergent power series.
Note again that $e^{a(x^i,\xi^i)}$ can be defined without assuming well-defined $\partial_t a(x^i,\xi^i)$.
Even without taking into account the detail property of the infinitesimal generator $\partial_t a(x^i,\xi^i)$, the exponentiability is realized by the boundedness of $a(x^i,\xi^i)$.
In this section, beginning with $a(x^i,\xi^i) = {\rm Log} ~ (U(x^i,\xi^i) + \kappa I)$, an algebraic module over a Banach algebra is defined.
The essential idea of presenting a useful algebraic structure is not to examine directly the set of $\partial_t a(x^i,\xi^i)$ but to focus on the set of $a(x^i,\xi^i)$ at first, and then the algebraic structure of infinitesimal generators is discovered in the next. 
Although $a(x^i,\xi^i) \in B(X)$ is trivially the infinitesimal generator, what is explained here is the structure of the set of pre-infinitesimal generators $\partial_t a(x^i,\xi^i)$.
 \vspace{2.5mm}\\

\begin{theorem}[Normed vector space]
Let $U_j(x^i,\xi^i)$ be evolution operators satisfying Eq.~(\ref{logex}), and ${\rm Log} ~ U_j(x^i,\xi^i)$ be well-defined for any $x^i,\xi^i \in [-L,L]$ and $j = 1,2, \cdots, n$.
${\rm Log} U_j(x^i,\xi^i)$ are assumed to commute with each other.
\[ \begin{array}{ll}
V_{Lg}(X) := \{ k {\rm Log} ~ U_j(x^i,\xi^i); ~ k \in {\bm C}, ~  x^i,\xi^i \in [-L,L] \} ~ \subset B(X)
  \subset G(X)  \end{array} \]
is a normed vector space\index{normed vector space} over the complex number field, where $B(X)$ denotes a set of all the bounded operators on $X$. \vspace{2.5mm}\\
\end{theorem}

{\bf Proof ~}
In case of $\kappa = 0$ the operator $a(x^i,\xi^i)$ is reduced to
\[ \begin{array}{ll}
 {\rm Log} ~ U(x^i,\xi^i) \in B(X).
\end{array} \]
The operator sum is calculated using the Dunford-Riesz integral
\begin{equation} \label{sum1} \begin{array}{ll}
  {\rm Log} ~ U(x^i,\eta^i) + {\rm Log} ~ U(\eta^i,\xi^i)   \vspace{2.5mm} \\
   = \frac{1}{2 \pi i} \int_{\Gamma}  {\rm Log} \lambda ~(\lambda I  - U(x^i,\eta^i))^{-1} d \lambda ~
    + \frac{1}{2 \pi i} \int_{\Gamma'}  {\rm Log} \lambda' ~(\lambda' I - U(\eta^i,\xi^i) )^{-1} d \lambda'  \vspace{2.5mm} \\
   = \frac{1}{(2 \pi i)^2}  \int_{\Gamma}  \int_{\Gamma'} ( {\rm Log} \lambda + {\rm Log} \lambda') 
     ~(\lambda I - U(x^i,\eta^i) )^{-1} ~(\lambda'I - U(\eta^i,\xi^i) )^{-1}
     ~ d \lambda'  d \lambda  \vspace{2.5mm} \\
   = \frac{1}{(2 \pi i)^2}  \int_{\Gamma}  \int_{\Gamma'} ( {\rm Log} \lambda  \lambda') 
     ~(\lambda I -  U(x^i,\eta^i))^{-1} ~(\lambda' I - U(\eta^i,\xi^i) )^{-1}
     ~ d \lambda'  d \lambda   \vspace{2.5mm} \\
   =   {\rm Log} ~ [ U(x^i,\eta^i) U(\eta^i,\xi^i)]
   =   {\rm Log} ~  U(x^i,\xi^i),
    \end{array} \end{equation}
then the sum closedness is clear.
Here $\Gamma'$ is assumed to be included in $\Gamma$, and this condition is not so restrictive in the present setting.
In a different situation, when $U(t,r)$ and $V(t,r)$ commute for the same $t$ and $r$, another kind of sum is calculated as
\begin{equation} \label{sum2} \begin{array}{ll}
  {\rm Log} ~ U_1(x^i,\eta^i) + {\rm Log} ~ U_2(x^i,\eta^i)  \vspace{2.5mm} \\
   = \frac{1}{2 \pi i} \int_{\Gamma}  {\rm Log} \lambda ~(\lambda I - U_1(x^i,\eta^i))^{-1} d \lambda ~
    + \frac{1}{2 \pi i} \int_{\Gamma'}  {\rm Log} \lambda' ~(\lambda' I - U_2(x^i,\eta^i))^{-1} d \lambda'  \vspace{2.5mm} \\
   = \frac{1}{(2 \pi i)^2}  \int_{\Gamma}  \int_{\Gamma'} ( {\rm Log} \lambda + {\rm Log} \lambda') 
     ~(\lambda I - U_1(x^i,\eta^i))^{-1} ~(\lambda' I - U_2(x^i,\eta^i))^{-1}
     ~ d \lambda'  d \lambda  \vspace{2.5mm} \\
   = \frac{1}{(2 \pi i)^2}  \int_{\Gamma}  \int_{\Gamma'} ( {\rm Log} \lambda  \lambda') 
     ~(\lambda I - U_1(x^i,\eta^i))^{-1} ~(\lambda' I - U_2(x^i,\eta^i))^{-1}
     ~ d \lambda'  d \lambda   \vspace{2.5mm} \\
   =   {\rm Log} ~ [ U_1(x^i,\eta^i) U_2(x^i,\eta^i)],
    \end{array} \end{equation}
where, for $W(x^i,\eta^i) = e^{ {\rm Log}U_1(x^i,\eta^i) + {\rm Log}U_2(x^i,\eta^i)}= U_1(x^i,\eta^i) U_2(x^i,\eta^i)$, the semigroup property is satisfied as
\[ \begin{array}{ll}
 W(x^i,\eta^i) W(\eta^i,\xi^i) = U_1(x^i,\eta^i) U_2(x^i,\eta^i) U_1(\eta^i,\xi^i) U_2(\eta^i,\xi^i)  = U_1(x^i,\eta^i) U_1(\eta^i,\xi^i) U_2(x^i,\eta^i)  U_2(\eta^i,\xi^i) = W(x^i,\xi^i), \vspace{2.5mm}\\
 W(\xi^i,\xi^i) = W(\xi^i,x^i) W(x^i,\xi^i) = I,
    \end{array} \]
and then the sum closedness is clear.
Although the logarithm function is inherently a multi-valued function, the uniqueness of sum operation is ensured by the single-valued property of the principal branch ``Log''.
Consequently, since the closedness for scalar product is obvious, where
\[
{\tilde W}(x^i,\eta^i) = e^{k_1 {\rm Log} ~ U_1(x^i,\eta^i)  + k_2 {\rm Log} ~ U_2(x^i,\eta^i)}
= e^{k_1 {\rm Log} ~ U_1(x^i,\eta^i) } e^{ k_2 {\rm Log} ~ U_2(x^i,\eta^i)} 
\]
holds the semigroup property in a similar way to $W(x^i,\eta^i) $. 
Consequently
$V_{Lg}(X)$
is a normed vector space over the complex number field.
In particular the zero operator ${\rm Log} I$ is included in $V_{Lg}(X)$.
Theorem 1 has been proved. \quad $\square$ \\

\begin{theorem}[$B(X)$-module]
Let $U_j(x^i,\xi^i)$ be evolution operators satisfying Eq.~(\ref{logex}) for any $x^i,\xi^i \in [-L,L]$ and $j = 1,2, \cdots, n$.
For a certain $K \in B(X)$, let a subset of $B(X)$ in which each element is assumed to commute with ${\rm Log} ~ (U_j(x^i,\xi^i)+K)$ be $B_{ab}(X)$.
${\rm Log} (U_j(x^i,\xi^i)+K)$ are assumed to commute with each other.
\[ \begin{array}{ll}
B_{Lg}(X) := \left\{ {\mathcal K}{\rm Log} ~ (U_j(x^i,\xi^i)+K); ~  {\mathcal K} \in B_{ab}(X),~ K \in B(X) ,~  x,\xi \in [-L,L] \right\} ~ \subset B(X) \subset G(X)
\end{array} \]
is a module over the Banach algebra. \vspace{2.5mm}\\
\end{theorem}

{\bf Proof ~}
It is worth generalizing the above normed vector space.
In this sense, utilizing a common operator $K \in B(X)$, components are changed to $ {\rm Log} ~ (U(x^i,\xi^i)+K)$.

The operator sum is calculated as
\begin{equation}  \begin{array}{ll}
  {\rm Log} ~ (U(x^i,\eta^i) + K) + {\rm Log} ~ (U(\eta^i,\xi^i)+ K)   \vspace{2.5mm} \\
   = \frac{1}{2 \pi i} \int_{\Gamma}  {\rm Log} \lambda ~(\lambda I - U(x^i,\eta^i)-K )^{-1} d \lambda ~
    + \frac{1}{2 \pi i} \int_{\Gamma'}  {\rm Log} \lambda' ~(\lambda' I - U(\eta^i,\xi^i)-K )^{-1} d \lambda'  \vspace{2.5mm} \\
   = \frac{1}{(2 \pi i)^2}  \int_{\Gamma}  \int_{\Gamma'} ( {\rm Log} \lambda + {\rm Log} \lambda') 
     ~(\lambda I - U(x^i,\eta^i)-K)^{-1} ~(\lambda' I - U(\eta^i,\xi^i)-K )^{-1}
     ~ d \lambda'  d \lambda  \vspace{2.5mm} \\
   = \frac{1}{(2 \pi i)^2}  \int_{\Gamma}  \int_{\Gamma'} ( {\rm Log} \lambda  \lambda') 
     ~(\lambda I - U(x^i,\eta^i)- K)^{-1} ~(\lambda' I - U(\eta^i,\xi^i)- K )^{-1}
     ~ d \lambda'  d \lambda   \vspace{2.5mm} \\
   =   {\rm Log} ~ [ (U(x^i,\eta^i) + K) (U(\eta^i,\xi^i) + K)]    \vspace{2.5mm} \\
   =   {\rm Log} ~ [ U(x^i,\xi^i) + K U(x^i,\eta^i) + K U(\eta^i,\xi^i) + K^2 ].
    \end{array} \end{equation}
After introducing a certain $K \in B(X)$ with sufficient large $\| K \|$, it is always possible to take integral path $\Gamma'$ to be included in $\Gamma$.
Since the part ``$K U(t,r) + K U(r,s) + K^2 I$'' is included in $B(X)$, the sum-closedness is clear. 
In a different situation, when $U_1(t,r)$ and $U_2(t,r)$ commute for the same $t$ and $r$, another kind of sum is calculated as
\begin{equation}  \begin{array}{ll}
  {\rm Log} ~ (U_1(x^i,\eta^i)+K_1) + {\rm Log} ~ (U_2(x^i,\eta^i)+K_2)   \vspace{2.5mm} \\
   = \frac{1}{2 \pi i} \int_{\Gamma}  {\rm Log} \lambda ~(\lambda I - U_1(x^i,\eta^i)- K_1 )^{-1} d \lambda ~
    + \frac{1}{2 \pi i} \int_{\Gamma'}  {\rm Log} \lambda' ~(\lambda' I - U_2(x^i,\eta^i)-K_2 )^{-1} d \lambda'  \vspace{2.5mm} \\
   = \frac{1}{(2 \pi i)^2}  \int_{\Gamma}  \int_{\Gamma'} ( {\rm Log} \lambda + {\rm Log} \lambda') 
     ~(\lambda I - U_1(x^i,\eta^i)-K_1 )^{-1} ~(\lambda'I - U_2(x^i,\eta^i)-K_2 )^{-1}
     ~ d \lambda'  d \lambda  \vspace{2.5mm} \\
   = \frac{1}{(2 \pi i)^2}  \int_{\Gamma}  \int_{\Gamma'} ( {\rm Log} \lambda  \lambda') 
     ~(\lambda I - U_1(x^i,\eta^i)-K_1 )^{-1} ~(\lambda' I - U_2(x^i,\eta^i)-K_2 )^{-1}
     ~ d \lambda'  d \lambda   \vspace{2.5mm} \\
  =   {\rm Log} ~ [ U_1(x^i,\eta^i) + K_1) (U_2(x^i,\eta^i) + K_2)]    \vspace{2.5mm} \\
  =   {\rm Log} ~ [ W(x^i,\eta^i) + K_2 U_1(x^i,\eta^i) + K_1 U_2(x^i,\eta^i) + K_1 K_2].
    \end{array} \end{equation}
Since the part ``$ K_2 U_1(x^i,\eta^i) + K_1 U_2(x^i,\eta^i) + K_1 K_2 $'' is included in $B(X)$, the sum-closedness is clear.

The product ${\mathcal K} {\rm Log} ~ (U_i(x^i,\xi^i) + K) \in B(X)$ is justified by the operator product equipped with $B(X)$. 
Since the closedness for operator product within $B(X)$ is obvious, where
\[
{\hat W}(x^i,\eta^i) = e^{{\mathcal K}_1 {\rm Log} ~ U_1(x^i,\eta^i)    + {\mathcal K}_2 {\rm Log} ~ U_2(x^i,\eta^i) }
= e^{{\mathcal K}_1 {\rm Log} ~ U_1(x^i,\eta^i)  } e^{ {\mathcal K}_2 {\rm Log} ~ U_2(x^i,\eta^i) } 
\]
holds the semigroup property, and using an identity operator $I \in B(X)$,
\[ \begin{array}{ll}
 {\mathcal K} {\rm Log} ~ (U(x^i,\eta^i) + K)  -  {\mathcal K} {\rm Log} ~ (U(x^i,\eta^i) )
 =  {\mathcal K}  {\rm Log} ~ [ I + K U(x^i,\eta^i) ^{-1}  ]
\end{array} \]
and therefore
\[ \begin{array}{ll}
  {\mathcal K}_1 {\rm Log} ~ (U_1(x^i,\eta^i) + K_1)    + {\mathcal K}_2 {\rm Log} ~ (U_2(x^i,\eta^i) +K_2)  
=  {\rm Log} ~ ( {\hat U} (x^i,\eta^i) + I ) + {\rm Log} ~ {\hat W}(x^i,\eta^i)  \vspace{2.5mm} \\
\quad = {\rm Log} ~[   {\hat W}(x^i,\eta^i) + {\hat U} (x^i,\eta^i)  {\hat W}(x^i,\eta^i) ]
\end{array} \]
are valid for $ {\hat U}(x^i,\eta^i) : = e^{ {\mathcal K}_1  {\rm Log} ~ [ I + K_1 U_1(x^i,\eta^i) ^{-1}  ] +  {\mathcal K}_2  {\rm Log} ~ [ I + K_2 U_2(x^i,\eta^i) ^{-1}  ] } - I$.
Consequently,
$B_{Lg}(X)$ 
is a module over a Banach algebra.
In particular a relation $V_{Lg}(X) \subset B_{Lg}(X)$ is satisfied.
The statement has been proved.  \quad $\square$ \\

The next corollary follows.

\begin{corollary}[$B(X)$-module for infinitesimal generators]
Let $U_j(x^i,\xi^i)$ be evolution operators satisfying Eq.~(\ref{logex}) for any $x^i,\xi^i \in [-L,L]$ and $j = 1,2, \cdots, n$.
For a certain $K \in B(X)$, let a subset of $B(X)$ in which each element is assumed to commute with ${\rm Log} ~ (U_j(x^i,\xi^i)+K)$ be $B_{ab}(X)$.
${\rm Log} (U_j(x^i,\xi^i)+K)$ are assumed to commute with each other.
\[ \begin{array}{ll}
G_{Lg}(X) := \left\{ {\mathcal K} \partial_{x^i} {\rm Log} ~ (U_j(x^i,\xi^i)+K); ~  {\mathcal K} \in B_{ab}(X),~ K \in B(X) ,~  x,\xi \in [-L,L] \right\} ~\subset G(X)
\end{array} \]
is a module over a Banach algebra. \vspace{2.5mm}\\
\end{corollary}

{\bf Proof ~}
According to the linearity of differential operator, the introduction of differential operator $\partial_{x_i} $ is true without any additional treatment.
It is sufficient to see that there exists a certain $\eta^i$ such that
\[ \begin{array}{ll}
 {\mathcal K} \partial_{x^i} {\rm Log} ~ (U_j(x^i,\xi^i)+K)
 =\partial_{x^i} \left[ \int_{\xi^i}^{x^i}  {\mathcal K} \partial_{\eta^i} {\rm Log} ~ (U_j(\eta^i,\xi^i)+K) d \eta^i  \right]       \vspace{1.5mm} \\
  =\partial_{x^i} \left[  {\mathcal K}  {\rm Log} ~ (U_j(x^i,\xi^i)+K) 
  - {\mathcal K}  {\rm Log} ~ (I +K) 
  -\left( \int_{\xi^i}^{x^i}   \partial_{\eta^i}  {\mathcal K} d \eta^i \right)  ~  {\rm Log} ~ (U_j( \eta^i,\xi^i)+K)  \right]   
\end{array} \] 
according to the mean value theorem (${\xi^i} \le  \eta^i \le {x^i}$).
 \quad $\square$ \\

The module over a Banach algebra is called $B(X)$-module\index{B(X)-module}.
For the structure of $B_{Lg}(X)$, a certain originally unbounded part can be classified to ${\rm Log} ~ (U(x^i,\xi^i)+K) \in B(X)$, and the rest part to ${\mathcal K} \in B_{ab}(X)$. 
Here the terminology ``originally unbounded'' is used, because some unbounded operators are reduced to bounded operators under the validity of the logarithmic representation.   

It is necessary to connect the concept of $B(X)$-module to the set of infinitesimal generators. 
Let us move on to the operator $\partial_{x^i} [ {\mathcal K}  {\rm Log} (U_j (x^i, \xi^i) + K) ]$, which is expected to be the pre-infinitesimal generator of $\exp [{\mathcal K}  {\rm Log} (U_j((x^i, \xi^i)) + K)]$.  
This property is surely true by the inclusion relation $B_{Lg}(X) \subset B(X)$.
It is also suggested by the inclusion relation $B_{Lg}(X) \subset B(X)$, operators $\partial_{x^i} {\rm Log} (U_j((x^i, \xi^i)) + K)]$ are the pre-infinitesimal generators if ${\rm Log} (U_j((x^i, \xi^i)) + K)] \in G_{Lg}(X)$ is satisfied.
Consequently the unbounded sum-perturbation for infinitesimal generators is seen by the sum closedness of $B(X)$-module. 
Note that it does not require the self-adjointness\index{self-adjointness} of the operator.

The pre-infinitesimal generator property is examined for products of operators in the next two theorems. \\

\begin{theorem}[Product perturbation for pre-infinitesimal generators] \label{bm-01}
For a certain $K \in B(X)$, let a subset of $B(X)$ in which each element is assumed to commute with ${\rm Log} ~ (U_i(x^i, \xi^i)+K)$ be $B_{ab}(X)$.
Let an operator denoted by 
\[
L  {\rm Log} (U(x^i, \xi^i) + K)
\]
be included in $B_{Lg}(X)$, where the evolution operator $U(t,s)$ on $X$ is generated by $A(t)$, $L$ is an element in $B_{ab}(X)$, and $K$ is an element in $B(X)$.
Let $K$ and $L$ be further assumed to be independent of $x^i$.
The product of pre-infinitesimal generators, which is represented by
\begin{equation}  \begin{array}{ll}
 L A(t) =   (I - K e^{-a(x^i, \xi^i)} )^{-1} \partial_{x^i} [ L  {\rm Log} (U(x^i, \xi^i) + K)],
\end{array} \end{equation}
is also the pre-infinitesimal generators in $X$, where $a(x^i, \xi^i) = {\rm Log} (U(x^i, \xi^i) + K)$.
\end{theorem}

{\bf Proof ~}
Since $L$ is independent of $x^i$,
\[
\partial_{x^i} [ L {\rm Log} (U(x^i, \xi^i) + K) ]
 =  L \partial_{x^i} [ {\rm Log} (U(x^i, \xi^i) + K)]
\]
is true.
The basic calculi using the $t$-independence of $K$ leads to the product of operator $L A(t)$.
It is well-defined by
\[  \begin{array}{ll}
    (I - K  e^{-a(x^i, \xi^i)})^{-1} \partial_{x^i} [ L {\rm Log} (U(x^i, \xi^i) + K)]  \vspace{1.5mm}  \\ 
\quad  = L (I -K  e^{-a(x^i, \xi^i)})^{-1}  \partial_{x^i} [ {\rm Log} (U(x^i, \xi^i) + K)]    \vspace{1.5mm}  \\
\quad  = L(I - K  e^{-a(x^i, \xi^i)})^{-1}   (U(x^i, \xi^i) + K)^{-1}  \partial_{x^i} U(x^i, \xi^i)   
 = L A(t)
\end{array} \]
under the commutation assumptions, where the relation $\partial_{x^i} U(x^i, \xi^i) = A (x^i)U(x^i, \xi^i)$ is applied.
Let $x^i,\xi^i \in [-L,+L]$ satisfy $\xi^i < x^i$. 
The pre-infinitesimal generator property of $L A(t)$ is confirmed by
\[  \begin{array}{ll}
 \left\| \int_{\xi^i}^{x^i}   (I - K e^{-a(x^i, \xi^i)})^{-1} \partial_{\eta^i} [ L {\rm Log} (U(\eta^i, \xi^i) + K)]   d\eta^i \right\| \vspace{1.5mm}  \\
  \le      \left\|     (I - K e^{-a(x^i, \xi^i)})^{-1} \int_{\xi^i}^{x^i}   \partial_{\eta^i} [ L {\rm Log} (U(\eta^i, \xi^i) + K)]   d\eta^i \right\| \vspace{1.5mm}  \\
  \le   {\displaystyle \sup_{\eta^i \in [\xi^i,x^i]} } \left\|   (I - K e^{-a(x^i, \xi^i)})^{-1}  \right\|  
 \left\|  \int_{\xi^i}^{x^i}   \partial_{\eta^i} [ L {\rm Log} (U(\eta^i, \xi^i) + K)]   d\eta^i \right\|  \vspace{1.5mm}  \\
 \le \Bigl( {\displaystyle \sup_{\eta^i \in [\xi^i,x^i]} } \left\|   (I - K e^{-a(x^i, \xi^i)})^{-1}  \right\|    \Bigr)
   \left\|     L  \right\|
   \left\|   {\rm Log} ((U(x^i, \xi^i)+ K) -  {\rm Log} (I + K) \right\|,
\end{array} \]
where a certain real number $\sigma \in [\xi^i,x^i]$ is determined by the mean value theorem. 
Consequently, due to the boundedness of $ \int_{\xi^i}^{x^i}   (I - K e^{-a(x^i, \xi^i)})^{-1} \partial_{x^i} [ L {\rm Log} (U(x^i, \xi^i) + K)]   d\eta^i $ on $X$, $(I - K e^{-a(x^i, \xi^i)})^{-1} \partial_{x^i} [ L {\rm Log} (U(x^i, \xi^i) + K)] $ is confirmed to be the pre-infinitesimal generator in $X$. \quad $\square$ \\

The operator $L$ can be regarded as a perturbation to the operators in $B_{Lg}(X)$.
This lemma shows the product-perturbation for the infinitesimal generators of $C^0$-semigroups under the commutation, although the perturbation has been studied mainly for the sum of operators.
It is remarkable that the self-adjointness of the operator is not required for this lemma.
For the details of conventional bounded sum-perturbation, and the perturbation theory for the self-adjoint operators, see Ref.~\cite{66kato}. \\

\begin{theorem}[Operator product] \label{bm-02}
For a certain $K(x^i) \in B(X)$, let a subset of $B(X)$ in which each element is assumed to commute with ${\rm Log} ~ (U(x^i, \xi^i)+K(x^i))$ be $B_{ab}(X)$.
Let an operator denoted by 
\[
L  {\rm Log} (U(x^i, \xi^i) + K(x^i))
\]
be included in $B_{Lg}(X)$, where the evolution operator $U(x^i, \xi^i)$ on $X$ is generated by $A(x^i)$, $L$ is an element in $B_{ab}(X)$, and $K(x^i)$ is an element in $B(X)$.
Let $L$ and $K(x^i)$ be $x^i$-independent and $x^i$-dependent, respectively.
The operators represented by
\begin{equation}  \begin{array}{ll}  \label{prorep}
{\mathcal L}(x^i)  \partial_{x^i} [ L  {\rm Log} (U(x^i, \xi^i) + K(x^i))]
\end{array} \end{equation}
is the pre-infinitesimal generators in $X$, if the operator ${\mathcal L}(\eta^i) \in B(X)$ is strongly continuous with respect to $\eta^i$ in the interval $[\xi^i,x^i]$. 
\end{theorem}

{\bf Proof ~}
Let $x^i, \xi^i \in [-L,+L]$ satisfy $\xi^i <x^i$. 
The pre-infinitesimal generator property is reduced to the possibility of applying the mean-value theorem.
\[  \begin{array}{ll}
 \left\| \int_{\xi^i}^{x^i}  {\mathcal L}(\eta^i)  \partial_{\eta^i} [   {\rm Log} (U(\eta^i, \xi^i) + K(\eta^i))]  d\eta^i \right\|   \vspace{1.5mm}  \\
  \le     {\mathcal L}(\sigma^i)    \left\| \int_{\xi^i}^{x^i} \partial_{\eta^i} [   {\rm Log} (U(\eta^i, \xi^i) + K(\eta^i))]  d\eta^i \right\|   \vspace{1.5mm}  \\
  \le   {\displaystyle \sup_{\sigma^i \in [\xi^i,x^i]} } \left\|  {\mathcal L}(\sigma^i)  \right\|   
   \left\| \int_{\xi^i}^{x^i} \partial_{\eta^i} [   {\rm Log} (U(\eta^i, \xi^i) + K(\eta^i))]  d\eta^i \right\|  \vspace{1.5mm}  \\
 \le    {\displaystyle \sup_{\sigma^i \in [\xi^i,x^i]} } \left\|  {\mathcal L}(\sigma^i)  \right\| 
   \left\|   {\rm Log} (U(x^i, \xi^i) + K(x^i)) -  {\rm Log} (I + K(\xi^i)) \right\|,
\end{array} \]
where a certain real number $\sigma^i \in [\xi^i,x^i]$ is determined by the mean value theorem. 
Consequently, ${\mathcal L}(x^i)  \partial_{x^i} [ L  {\rm Log} (U(x^i, \xi^i) + K(x^i))]$ is confirmed to be the pre-infinitesimal generator in $X$. \quad $\square$ \\

Equation (\ref{prorep}) provides one standard form for the representation of operator products in the sense of logarithmic representation. 
Consequently $B(X)$-module is associated with the pre-infinitesimal generator.\\

\subsection{Formulation of rotation group}
The application example of $B(X)$-module is provided.
The concept of $B(X)$-module is general enough to provide a foundation of the conventional bounded formulation of Lie algebra\index{Lie algebra}s (for a textbook, see \cite{73sagle}).
In other words, by means of $B(X)$-module, the intersection of the Banach algebra (including only bounded operators) and the extracted bounded part of the Lie algebra (generally including unbounded operators) is shown.
More precisely, using $B(X)$-module, the bounded part is extracted from unbounded angular momentum operators.
The extracted bounded parts are utilized to formulate the rotation group with incorporating the unboundedness of angular momentum algebra.

The mathematical foundation of rotation group\index{rotation group} is demonstrated~\cite{19iwata}.
Although the evolution parameter in this paper is denoted by $t,s \in [-T,+T]$, it is more likely to be denoted by $\theta, \sigma \in [-\Theta, +\Theta]$, because the evolution parameter in the present case means the rotation angle.
The rotation group is generated by the angular momentum operator (for textbooks, see Refs.~\cite{50weyl,57yamanouchi}).
The angular momentum operator\index{angular momentum operator} includes a differential operator, as represented by
\begin{equation} \begin{array}{ll}
{\mathcal L} = -i \hbar ({\bf r} \times \nabla),
\end{array} \end{equation}
where $\hbar$ is a real constant called the Dirac constant. 
The appearance of differential operator $\nabla$ in the representation of ${\mathcal L}$ is essential.
The operator $\nabla$ is an unbounded operator for example in a Hilbert space $L^2({\bf R}^3)$, while it must be treated as a bounded operator in terms of establishing an algebraic ring structure.
Furthermore, the operator boundedness is also indispensable for some important formulae such as the Baker-Campbell-Hausdorff formula\index{Baker-Campbell-Hausdorff formula} and the Zassenhaus formula\index{Zassenhaus formula} to be valid.
In general, the exponential of unbounded operators cannot be represented by the power series expansion (cf. the Yosida approximation in a typical proof of the Hille-Yosida theorem; e.g., see Ref.~\cite{79tanabe}).

Let ${\mathbb R}^3$ be the three-dimensional spatial coordinate spanned by the standard orthogonal axes, $x$, $y$ and $z$.
The angular momentum operator ${\mathcal L}$ is considered in $L^2({\mathbb R}^3)$.
The angular momentum operator
\[ \begin{array}{ll}
{\mathcal L} =  ({\mathcal L}_x, {\mathcal L}_y,{\mathcal L}_z)  \vspace{2.5mm} \\
\end{array} \] 
consists of $x$, $y$, and $z$ components
\[  \begin{array}{ll}
\quad {\mathcal L}_x  = -i \hbar ( y \partial_{z} -  z \partial_{y}),    \vspace{2.5mm} \\
\quad {\mathcal L}_y  = -i \hbar ( z \partial_{x} -  x \partial_{z}),   \vspace{2.5mm} \\
\quad {\mathcal L}_z  = -i \hbar ( x \partial_{y} -  y \partial_{x}), 
\end{array} \]
respectively.
The commutation relations
\begin{equation} \label{eqcom}  \begin{array}{ll}
[{\mathcal L}_x, {\mathcal L}_y] = i \hbar {\mathcal L}_z, \quad
[{\mathcal L}_y, {\mathcal L}_z] = i \hbar {\mathcal L}_x, \quad
[{\mathcal L}_z, {\mathcal L}_x] = i \hbar {\mathcal L}_y
\end{array} \end{equation}
are true, where $[{\mathcal L}_i,{\mathcal L}_j]:={\mathcal L}_i {\mathcal L}_j-{\mathcal L}_j {\mathcal L}_i$ denotes a commutator product\index{commutator product} ($i,j =x,y,z$).
The commutation of angular momentum operators arises from the commutation relations of the canonical quantization
\begin{equation} \label{omt}  \begin{array}{ll}
[x, p_x] = [y, p_y] = [z, p_z] = i \hbar,  \vspace{1.5mm}  \\
\left[y, p_x \right] = \left[y, p_z\right] = \left[z, p_x\right] = \left[z, p_y\right] = \left[x, p_y\right] = \left[x, p_z\right]  = 0.
\end{array} \end{equation}
Indeed, the momentum operator $p=(p_x,p_y,p_z)$ is represented by $p=-i \hbar ( \partial_{x},  \partial_{y},  \partial_{z})$ in quantum mechanics.
It is remarkable that the commutation is always true for the newtonian mechanics;  i.e., $[x, p_x] = [y, p_y] = [z, p_z] =0$ is true in addition to $\left[y, p_x \right] = \left[y, p_z\right] = \left[z, p_x\right] = \left[z, p_y\right] = \left[x, p_y\right] = \left[x, p_z\right]  = 0$.

Let a set of all bounded operators on $L^2({\mathbb R}^3)$ be denoted by $B(L^2({\mathbb R}^3))$.
A set of operators $\{ {\mathcal L}_k; ~ k = x,y,z \}$ or  $\{ i {\mathcal L}_k/\hbar; ~ k = x,y,z \}$ with the commutation relation (\ref{eqcom}) is regarded as the Lie algebra. 
In particular $\{ {\hat \alpha} {\mathcal L}_x + {\hat \beta} {\mathcal L}_y + {\hat \gamma} {\mathcal L}_z; ~ {\hat \alpha}, {\hat \beta}, {\hat \gamma} \in {\mathbb C} \}$ forms a vector space over the complex number field, while $\{ \alpha( i {\mathcal L}_x/\hbar) + \beta (i {\mathcal L}_y/\hbar) + \gamma (i {\mathcal L}_z/\hbar); ~ \alpha, \beta, \gamma \in {\mathbb R} \}$ is a vector space over the real number field.
It is possible to associate the real numbers $\alpha$, $\beta$, and $\gamma$ with the Euler angles (for example, see Ref.~\cite{82landau}).
The second term of the right hand side of
\begin{equation} \begin{array}{ll} \label{intem}
  ( r_i \partial_{j} ) ( r_k \partial_{l} ) 
  =  r_i (\partial_{j}   r_k) \partial_{l} +  r_i   r_k \partial_{j} \partial_{l} 
\end{array} \end{equation}
disappears as far as the commutator product $[{\mathcal L}_i,{\mathcal L}_j]$ is concerned, where $r_i$ is equal to $i$, and $i,j,k,l =x$, $y$, or $z$ satisfy $i \ne j$ and $k \ne l$.
This fact is a key to justify the algebraic ring structure of $\{ {\mathcal L}_k ; ~ k = x,y,z \}$.
On the other hand, although ${\mathcal L}_k$ is assumed to be bounded on $L^2({\mathbb R}^3)$ in the typical treatment of the Lie algebra, it is not the case for the angular momentum algebra because of the appearance of differential operators in their definitions. 
From a geometric point of view, the range space $R({\mathcal L_k}) \subset L^2({\mathbb R}^3)$ strictly includes the domain space $D({\mathcal L_k})$; i.e., there is no guarantee for any $u \in L^2({\mathbb R}^3)$ and a certain positive $M \in {\mathbb R}$ to satisfy $ \| {\mathcal L_k} u \|_{ L^2({\mathbb R}^3)} \le M \| u \|_{ L^2({\mathbb R}^3)}$.
In order to establish $\{ i {\mathcal L}_k/\hbar; ~ k = x,y,z \}$ as the Lie algebra, it is necessary to show 
\[
\pm  i {\mathcal L}_k/\hbar =  \pm (r_i \partial_{r_j} - r_j \partial_{r_i})
\]
as an infinitesimal generator in $L_{r}^2({\mathbb R}^3)$, where $i,j,k =x,y,z$ satisfies $i \ne j \ne k$. 
As for the angular momentum operator, the $t$-independent assumption for operators $K$ and $L$ in Theorem~\ref{bm-01} and Theorem~\ref{bm-02} is satisfied.
Note that $t$-independence assumes a kind of commutation relation. \\

\begin{theorem}[Unbounded formulation of rotation group]
Let $r_i$ be either $x$, $y$, or $z$.
For $i \ne j$, an operator $\pm r_i \partial_{r_j}$ with its domain space $H_{r}^1({\mathbb R}^3)$ is an infinitesimal generator in $L_{r}^2({\mathbb R}^3)$.
Consequently, the angular momentum operators 
\[  \begin{array}{ll}
\quad \pm i {\mathcal L}_x/\hbar  = \pm ( y \partial_{z} -  z \partial_{y}),    \vspace{2.5mm} \\
\quad \pm i {\mathcal L}_y/\hbar  = \pm  ( z \partial_{x} -  x \partial_{z}),   \vspace{2.5mm} \\
\quad \pm i {\mathcal L}_z/\hbar  = \pm ( x \partial_{y} -  y \partial_{x})
\end{array} \]
are infinitesimal generators in $L_{r}^2({\mathbb R}^3)$.
\end{theorem}

{\bf Proof ~} 
The proof consists of three steps. \\

{\bf [1st step: $\partial_{r_k}$ as an infinitesimal generator] ~}
\begin{lemma} \label{lem-r1}
For $r_k$ equal to $x$, $y$, or $z$, an operator $\partial_{r_k} $ with its domain $H^1({\mathbb R}^3)$ is an infinitesimal generator in $L^2({\mathbb R}^3)$.
\end{lemma}

 {\bf [Proof of Lemma \ref{lem-r1}].} 
The operator $\partial_x$ is known as the infinitesimal generator of the first order hyperbolic type partial differential equations.
For a complex number $\lambda$ satisfying ${\rm Re} \lambda > 0$, let us consider a differential equation
\begin{equation} \label{de1} \begin{array}{ll}
 \partial_x  u(x) =  \lambda u(x) - f(x) 
\end{array} \end{equation}
in $L^2({\mathbb R})$, and
\[ \begin{array}{ll}
 u(x) = - \int_x^{\infty}  \exp [\lambda (x-\xi)] f(\xi) d \xi
\end{array} \]
satisfies the equation.
According to the Schwarz inequality,
\[ \begin{array}{ll}
 \int_{-\infty}^{+\infty} | u(x) |^2 dx
  = \int_{-\infty}^{+\infty} | \int_x^{\infty}  \exp [\lambda (x-\xi)] f(\xi) d \xi |^2 dx \vspace{1.5mm}  \\
\qquad  \le  \int_{-\infty}^{+\infty} \left\{ \int_x^{\infty} \exp \left[ \frac{({\rm Re} \lambda) (x-\xi)}{2} \right]  \exp \left[ \frac{({\rm Re} \lambda) (x-\xi)}{2} \right] |f(\xi)| d \xi \right\}^2 dx  \vspace{1.5mm}  \\
\qquad  \le  \int_{-\infty}^{+\infty} \int_x^{\infty} \exp \left[({\rm Re} \lambda) (x-\xi) \right] d \xi ~
  \int_x^{\infty} \exp \left[({\rm Re} \lambda) (x-\xi) \right] |f(\xi)|^2 d \xi ~ dx
\end{array} \]
is obtained, because $|e^{\lambda/2}|^2 = |e^{{\rm Re}\lambda/2}|^2 ~ |e^{i{\rm Im}\lambda/2}|^2 \le e^{{\rm Re}\lambda}$ is valid if ${\rm Re} \lambda > 0$.
Here the equality
\[ \begin{array}{ll}
 \int_x^{\infty} \exp \left[({\rm Re} \lambda) (x-\xi) \right] d \xi
 =  \int_0^{\infty} \exp \left[(-{\rm Re} \lambda) \xi \right] d \xi 
 =  \frac{1}{{\rm Re} \lambda}
\end{array} \]
is positive valued if ${\rm Re} \lambda > 0$.
Its application leads to
\[ \begin{array}{ll}
 \int_{-\infty}^{+\infty} | u(x) |^2 dx
 \le \frac{1}{{\rm Re} \lambda}  \int_{-\infty}^{+\infty}  \int_x^{\infty} \exp \left[({\rm Re} \lambda) (x-\xi) \right] |f(\xi)|^2 d \xi~ dx  \vspace{1.5mm}\\
\qquad \le \frac{1}{{\rm Re} \lambda}  \int_{-\infty}^{+\infty}   \int_{-\infty}^{\xi} \exp \left[({\rm Re} \lambda) (x-\xi) \right] dx  |f(\xi)|^2 d \xi.
\end{array} \]
Further application of the equality
\[ \begin{array}{ll}
 \int_{-\infty}^{\xi} \exp \left[({\rm Re} \lambda) (x-\xi) \right] dx
 =  \int_{-\infty}^0 \exp \left[({\rm Re} \lambda) x \right] dx
 =  \frac{1}{{\rm Re} \lambda}
\end{array} \]
results in
\[ \begin{array}{ll}
 \int_{-\infty}^{+\infty} | u(x) |^2 dx
 \le \frac{1}{{\rm Re} \lambda^2}  \int_{-\infty}^{+\infty}  |f(\xi)|^2  d \xi,
\end{array} \]
and therefore
\[ \begin{array}{ll}
\| (\lambda I - \partial_x)^{-1} f \|_{L^2({\mathbb R})} \le  \frac{1}{{\rm Re} \lambda^2} \| f \|_{L^2({\mathbb R})}.
\end{array} \]
That is, for ${\rm Re} \lambda > 0$,
\[ \begin{array}{ll}
\| (\lambda I - \partial_x)^{-1}  \| \le  \frac{1}{{\rm Re} \lambda}
\end{array} \]
is valid. 
The surjective property of $(\lambda I - \partial_k)$ is seen by the unique existence of solution $u(x) \in L^2({\mathbb R})$ for the initial value problem of Eq.~(\ref{de1}).

A semigroup is generated by taking a subset of the complex plane as
\[ \begin{array}{ll}
\Omega = \{ \lambda \in {\mathbb C}; ~ \lambda = \overline{\lambda} \}
\end{array} \]
where $\Omega$ is included in the resolvent set of $\partial_x$.
For $\lambda \in \Omega$, $(\lambda I - \partial_x)^{-1}$ exists, and
\[ \begin{array}{ll}
\| (\lambda I - \partial_x)^{-n}  \| \le  \frac{1}{({\rm Re} \lambda)^n}
\end{array} \]
is obtained.
Consequently, according to the Lumer-Phillips theorem \cite{61lumer,52phillips} for the generation of quasi contraction semigroup, $\partial_x$ with the domain space $H^1({\mathbb R})$ is confirmed to be an infinitesimal generator in $L^2({\mathbb R})$.
The similar argument is valid to $\partial_y$ and  $\partial_z$.
By considering $(x,y,z) \in {\mathbb R}^3$,  $\partial_k$ with $k=x,y,z$ are the infinitesimal generators in $L^2({\mathbb R}^3)$. \quad 
 {\bf [Q.E.D.: Lemma \ref{lem-r1}]}    \vspace{2.5mm}\\

{\bf [2nd step: $i r_k I$ as an infinitesimal generator] ~}
\begin{lemma} \label{lem-r2}
Let $r_k$ be either $x$, $y$, or $z$.
Let $I$ be the identity operator of $L^2({\mathbb R}^3)$.
An operator $i r_k I$ is an infinitesimal generator in $L^2({\mathbb R}^3)$.
\end{lemma}

 {\bf [Proof of Lemma \ref{lem-r2}].}  
For any $w \in {\mathbb C}$, it is possible to define the exponential function by the convergent power series:
$e^{w} = {\displaystyle \sum_{j=0}^{\infty}} (w)^j/j!$, so that 
\[  \begin{array}{ll}
 e^{i t r_k I} ={\displaystyle \sum_{j=0}^{\infty}} \frac{1}{j!} (i t r_k I)^j 
  \end{array} \]
is well-defined for $t, r_k \in {\mathbb R}$.
This fact is ensured by the boundedness of the identity operator $I$, although $r_k I$ and $i r_k I$ are not bounded operators in $L^2({\mathbb R})$ if the standard $L^2$-norm is equipped.
It is sufficient for $i r_k I$ to be the pre-infinitesimal generator.

For an arbitrary $r_k \in {\mathbb R}$, an operator $i r_k I$ with its domain $L^2({\mathbb R})$ is the infinitesimal generator in $L^2({\mathbb R})$; 
indeed, the spectral set is on the imaginary axis of the complex plane, and the unitary operator is generated as
\[   \begin{array}{ll}
\int |(e^{it r_k I} u)|^2 ~ d r_k 
= \int |u|^2 ~ d r_k.
\end{array}  \]
Consequently, the operator $i r_k I$ is treated as an infinitesimal generator in $L^2({\mathbb R})$ and therefore in $L^2({\mathbb R}^3)$. \quad 
 {\bf [Lemma \ref{lem-r2}]}   \\

{\bf [3rd step: $\pm ( r_i \partial_{r_j} - r_j \partial_{r_i})$ as an infinitesimal generator] ~}
Let $i \ne j$ be satisfied for $i,j =x,y,z$.
Since $e^{\pm(t-s) \partial_{r_j}}$ is well-defined (cf. Lemma~\ref{lem-r1}) with the domain space $H_{r}^1({\mathbb R}^3)$, its logarithmic representation is obtained by
\[ \begin{array}{ll}
\pm \partial_{r_j} =  ( I + \kappa e^{\pm (s-t) \partial_{r_j}}) \partial_t {\rm Log} (e^{ \pm (t-s) \partial_{r_j}} + \kappa I),  
\end{array} \]
where $\kappa \ne 0$ is a certain complex number. 
According to Theorem~\ref{bm-02}, the product between $i r_i I$ and $ \pm \partial_{r_j}$ is represented by
\[
\pm i r_i \partial_{r_j} 
= i r_i   ( I - \kappa e^{ \pm (s-t) \partial_{r_j}}) \partial_t [ {\rm Log} (e^{ \pm (t-s) \partial_{r_j}} + \kappa I) ].
\]
Using the commutation and $t$-independence of $r_i I$, it leads to the logarithmic representation
\[
\pm  r_i \partial_{r_j} 
 =  ( I + \kappa e^{ \pm (s-t) \partial_{r_j}}) \partial_t [  r_i {\rm Log} (e^{ \pm (t-s) \partial_{r_j}} + \kappa I)]
\]
without the loss of generality.
The domain space of $r_i I$ is equal to $L^2({\mathbb R}^3)$, as $e^{ i r_i I}$ is represented by the convergent power series in $L_{r}^2({\mathbb R}^3)$.
The half plane $\{ \lambda \in {\mathbb C}; ~ {\rm Re} \lambda > 0 \}$ is included in the resolvent set of $  \pm r_i \partial_{r_j} $. 
Consequently, for $t, r_i \in {\mathbb R}$, the existence of $e^{\pm t r_i \partial_{r_j}}$ directly follows from the confirmed existence of $e^{\pm t \partial_{r_j}}$.
Being equipped with the domain space $H_{r}^1({\mathbb R}^3)$,  $\pm r_i \partial_{r_j}$ is the infinitesimal generator in $L_{r}^2({\mathbb R}^3)$. 

The pre-infinitesimal generator property of sum is also understood by the $B(X)$-module property.
The sum  between $r_i \partial_{r_j} $ and $- r_j \partial_{r_i} $ is represented by
\begin{equation} \begin{array}{ll} \label{sumre}
 ( I + \kappa e^{+ (s-t) \partial_{r_j}})  \partial_t  [ r_i ~ {\rm Log} (e^{+ (t-s) \partial_{r_j}} + \kappa I)]
 - 
 ( I + \kappa e^{- (s-t) \partial_{r_i}}) \partial_t  [ r_j ~ {\rm Log} (e^{- (t-s) \partial_{r_i}} + \kappa I)]  \vspace{1.5mm} \\
 =
  ( I + \kappa e^{+ (s-t) ~ \partial_{r_j}})  
   \partial_t  [ r_i ~ {\rm Log} (e^{+ (t-s) \partial_{r_j}} + \kappa I)  -     r_j ~ {\rm Log} (e^{- (t-s) \partial_{r_i}} + \kappa I)]
  \\
\quad - 
 ( \kappa e^{+ (s-t) ~ \partial_{r_j}}  -  \kappa e^{- (s-t) \partial_{r_i}}) \partial_t  [ r_j ~ {\rm Log} (e^{- (t-s) \partial_{r_i}} + \kappa I)]  \vspace{1.5mm} \\
 =
  ( I + \kappa e^{+ (s-t) ~ \partial_{r_j}}) 
   \partial_t   r_i [  {\rm Log} (e^{+ (t-s) \partial_{r_j}} + \kappa I)  -~ {\rm Log} (e^{- (t-s) \partial_{r_i}} + \kappa I) ]  \\
\quad -
  ( I + \kappa e^{+ (s-t) ~ \partial_{r_j}})  
   \partial_t  [ (r_i   -     r_j )~ {\rm Log} (e^{- (t-s) \partial_{r_i}} + \kappa I)]  \\
\quad - 
 ( \kappa e^{+ (s-t) ~ \partial_{r_j}}  -  \kappa e^{- (s-t) \partial_{r_i}}) \partial_t  [ r_j ~ {\rm Log} (e^{- (t-s) \partial_{r_i}} + \kappa I)],
\end{array} \end{equation}
where all the three terms in the right hand side are of the form
\[  \begin{array}{ll}
{\mathcal L}(t)  \partial_t [ r_i  {\rm Log} (U(t, s) + K(t))]
\end{array} \]
whose pre-infinitesimal generator properties are proved similarly to Lemma 2.
In particular, the first term in the right hand side can be reduced to the above form with $L=1$ and $t$-dependent $K(t)$, the parts corresponding to ${\mathcal L}(t)$ are strongly continuous, and $r_i$ is independent of $t$.  
After having an integral of Eq.~(\ref{sumre}) in terms of $t$, each term is regraded as a bounded operator on $L_{r}^2({\mathbb R}^3)$.
Consequently, for $i \ne j$, the application of Lemma 2 leads to the fact that
\[  \begin{array}{ll}
\pm ( r_i \partial_{r_j} - r_j \partial_{r_i})
\end{array} \]
with its domain space $H_{r}^1({\mathbb R}^3)$ is the infinitesimal generator in $L_{r}^2({\mathbb R}^3)$.  \quad $\square$ \\

\begin{corollary}[Collective renormalization]
For $t, s \in [-T, +T]$, let $V_k(t,s)$ with $k=x,y,z$ in $L_{r}^2({\mathbb R}^3)$ be generated by $i {\mathcal L}_k/\hbar$. 
For a certain complex constant $\kappa \ne 0$, the angular momentum operator $\pm i {\mathcal L}_k/\hbar $ with $k=x,y,z$ is represented by the logarithm
\begin{equation} \label{replk}
\pm i {\mathcal L}_k/\hbar = \pm ( I + \kappa V_k(s,t) ) \partial_t  [  {\rm Log} (V_k(t,s) + \kappa I)],
 \end{equation}
and the corresponding evolution operator is expanded by the convergent power series
\begin{equation} \begin{array}{ll} \label{replke}
V_k(t,s) = e^{ {\rm Log} (V_k(t,s) + \kappa I)} - \kappa I  
= {\displaystyle \sum_{n=0}^{\infty}} \frac{1}{n!}  ( {\rm Log} (V_k(t,s) + \kappa I) )^n - \kappa I  
=  (1-\kappa)I + {\displaystyle \sum_{n=1}^{\infty}} \frac{1}{n!}  ( {\rm Log} (V_k(t,s) + \kappa I) )^n 
\end{array} \end{equation}
where ${\rm Log} (V_k(t,s) + \kappa I)$ is bounded on $L_{r}^2({\mathbb R}^3)$, although ${\mathcal L}_k$ is unbounded in $L_{r}^2({\mathbb R}^3)$.
\end{corollary}

{\bf Proof ~}
The group $V_k(t,s)$ with $k=x,y,z$ is generated by the infinitesimal generator $i {\mathcal L}_k/\hbar$ in $L_{r}^2({\mathbb R}^3)$.
This fact leads to the logarithmic representation 
\begin{equation}
i {\mathcal L}_k/\hbar =  ( I + \kappa V_k(s,t) ) \partial_t  [  {\rm Log} (V_k(t,s) + \kappa I)],
 \end{equation}
where $\kappa \ne 0$ is a certain complex constant.
The relation $V_k(t,s) + \kappa I =  e^{ {\rm Log} (V_k(t,s) + \kappa I )}$ admits the power series expansion of $V_k(t,s)$.\quad  $\square$ \\

Equation~(\ref{replk}) shows a convergent power-series representation for the rotation group.
Let us call the representation shown in Eq.~(\ref{replk}) the collective renormalization\index{collective renormalization} (cf. renormalized evolution equation in Corollary~\ref{transform}), in which a detailed  degree of freedom $r_i \partial_{r_j}$ is switched to a collective degree of freedom ${\mathcal L}_k$.
According to the collective renormalization, the evolution problem is studied by beginning with the bounded evolution operator $V_k(s,t)$ and the related bounded infinitesimal generator ${\rm Log} (V_k(t,s) + \kappa I )$. 
In a more mathematical sense, the collective renormalization plays a role of simplifying the representation.
Equation (\ref{replke}) ensures the validity of convergent power series expansions used in operator algebras even if they include unbounded operators.


\section{Concluding remarks}
\subsection{Template of solvable nonlinear equations}
The utility of the logarithmic representation is found in a formal discussion.
The derivative of the logarithmic representation is formally represented by
\begin{equation} \label{leib} \begin{array}{ll}
\partial_t {\rm Log} v = v'  v^{-1},
\end{array} \end{equation}
where $v$ is a function of $t$, and the notation $'$ denotes the differentiation along the $t$-direction.
Since the logarithmic derivative $\partial_t {\rm Log} v$ corresponds to the infinitesimal generator if $v$ is the evolution operator, this equality shows the relation between the infinitesimal generator $\partial_t {\rm Log} v$ and the evolution operator $v$.
Let $v$ be a known function (possibly a solution of linear equation), and $u$ be an unknown function of $t$ ($u v^{-1}$ be a solution of another equation, and  of possibly a nonlinear equatoin).
The Leibnitz rule reads
\begin{equation} \label{leib2} \begin{array}{ll}
(u v^{-1})' =  [ v' v^{-1} -  u'  u^{-1} ] (u v^{-1}).
\end{array} \end{equation}
Both $v' v^{-1}$ and $u' u^{-1}$ are regarded as the logarithmic derivative for $t$-direction.
The change of the evolution direction simply requires to fix $u = \partial_x v$, and then $u v^{-1}$ is regarded as the logarithmic derivative for $x$-direction. 
As seen in the case of Cole-Hopf transform, Eq.~(\ref{leib2}) being equivalent to the Burger's equation in case of the Cole-Hopf transform provides one abstract template for nonlinear evolution equations, which can be analyzed as the linear problem.  
If the Cole-Hopf transform ($\psi = v'v^{-1}$) is combined with the Miura transform $w = \psi' + \psi^2$, the higher order version of Eq.~(\ref{leib})
\begin{equation}
  w = v''  v^{-1}
\end{equation}
is obtained~\cite{20iwata}.
In this way the logarithmic representation provides templates of solvable nonlinear equations, which can be reduced to linear equations.

\subsection{Related topics}
As for the applicability of the theory, the conditions to obtain the logarithmic representation (conditions shown in Sec.~\ref{tp-group}) are not so restrictive; indeed, they can be satisfied by $C_0$-semigroups generated by $x^i$-independent infinitesimal generators. 
The most restrictive condition to obtain the logarithmic representation is the commutation between $K(x^i)$ and $U(x^i,\xi^i)$.
Such a commutation is trivially satisfied by $x^i$-independent $K(x^i) = K$, and also satisfied when the variable $x^i$ is separable  (i.e., for an integrable function $g(x^i)$, $K(x^i) = g(x^i) K$).
In this sense the operator specified in Theorem~\ref{thm-rel} corresponds to a moderate generalization of $x^i$-independent infinitesimal generators.
The summary is demonstrated along with the related topics.  \vspace{3mm} \\
\underline{\bf Time reversal symmetry} ~
Let the existence of negative time evolution be a kind of time reversal symmetry.
Note that this kind of symmetry is true for linear wave equations, but false for linear heat equations.
The logarithmic representation of infinitesimal generators has been originally obtained for the invertible evolution operators, and it is generalized to non-invertible evolution operators.
Under the validity of boundedness of $U(x^i,\xi^i)$ on $X$, the removal of invertible criterion is essentially realized by the introduction of nonzero $\kappa \in {\mathbb C}$.
On the other hand, the indispensable conditions for obtaining this kind of logarithmic representations are the boundedness of the spectral set of $U(x^i,\xi^i)$ and the commutation assumption, where the bounded interval $-L \le x^i,\xi^i \le L$ is also necessary.
Consequently the time-reversal symmetry is recovered for the regularized evolution operator if $x^i$ is equal to $x^0$.
In the same way, a similar concept to spatial reversal symmetry being defined by the negative evolution can be recovered and violated by taking $x^i \ne x^0$. 
\vspace{3mm} \\
\underline{\bf Regularity} ~
The recovery of local time-reversal symmetry is associated with the regularity of the solution.
The concept of regularized trajectory, whose regularity is similar to that of the analytic semigroups (for a textbook, see Ref.~\cite{79tanabe}) at the least, is true for regularized evolution operators.
\vspace{3mm} \\
\underline{\bf Nonlinearity} ~
For obtaining the logarithmic representation, the operator $U(x^i,\xi^i)$ can be either linear or nonlinear semigroup.
The nonlinearity of semigroup can appear simply by altering the evolution direction under a suitable identification between the infinitesimal generator and the evolution operator. 
In particular the relation between evolution operator and its infinitesimal generator is essentially similar to the Cole-Hopf transform. \vspace{3mm} \\
\underline{\bf The self-adjointness} ~
The results obtained for a $B(X)$-module does not require the self-adjointness of the operator, so that it opens up a way to have a full-complex analysis (neither real nor pure-imaginary analysis) for a class of unbounded operators in association with the operator algebra.
It is worth noting that the obtained algebraic structure corresponds to a generalization of ``perturbation theory for semigroups of operators~\cite{66kato}''. \vspace{3mm} \\
\underline{\bf Discrete property} ~
For example, in case of two-dimensional space-time distribution, let the $C^0$-semigroup for $x^0$ direction exist for a Cauchy problem:
\[ \begin{array}{ll}
\partial_{x^0} U(x^0, \xi^0)~u_0  = K(x^0) U(x^0,\xi^0) ~u_0,  \vspace{2.5mm} \\
 \partial_{x^0} a(x^0,\xi^0) = \partial_{x^0} {\rm Log}  (U(x^0,\xi^0) + \kappa I),
\end{array} \]
in $X_0 := L^2(-L,L)$, where a certain complex number $\kappa$ is taken from the resolvent set of $U(x^i,\xi^i)$.
Furthermore let the same equation possible to be written as
\[ \begin{array}{ll}
\partial_{x^1} V(x^1, \xi^1)~v_0  = {\mathcal K}(x^1) V(x^1,\xi^1) ~v_0,  \vspace{2.5mm} \\
 \partial_{x^1} \alpha(x^1,\xi^1) = \partial_{x^1} {\rm Log}  (V(x^1,\xi^1) + \kappa I),
\end{array} \]
in $X_1:= L^2(-T,T)$.
In this situation, using $\alpha(x^1,\xi^1)$ instead of $a(x^0,\xi^0)$, the corresponding dynamical system holds a discrete trajectory in $X_0$ (for a illustration, see Fig.2 of Ref.~\cite{18iwata-2}).
Indeed, the trajectory is $L^2$-function with respect to $x^0$, and $C^0$ function with respect to $x^1$. 
That is, the relativistic treatment naturally leads to the discrete evolution (for a theory including the discrete evolution, see the variational method of abstract evolution equation \cite{61lions,72lions}).
The discrete evolution to the $t$-direction ($x^0$-direction), which can be obtained by altering the evolution direction, is expected to be useful to analyze the stochastic differential equations within the semigroup theory of operators.

\newpage
\section*{Data availability statement}
All datasets generated for this study are included in the article/supplementary material.

\section*{Conflict of interest}
The author declares that the research was conducted in the absence of any commercial or financial relationships that could be constructed as a potential conflict of interest. 

\section*{Acknowledgement}
The author is grateful to Prof. Emeritus Hiroki Tanabe for valuable comments.
This work was partially supported by JSPS KAKENHI Grant No. 17K05440.
Comments and suggested sentences from referees are appreciated.
\vspace{36mm} \\

\newpage

\printindex


\newpage
\begin{thebibliography}{1}
 \bibitem{02arendt}
W. Arendt,
Semigroups and evolution equations: functional calculus, regularity and kernel estimates, 
Handbook of Differential Equations {\bf 1}, 2002. \\

 \bibitem{16arendt}
W. Arendt, D. Dier, and S. Fackler, 
J. L. Lions's Problem on Maximal Regularity,
arXiv:1612.03676, 2016. \\

\bibitem{15bateman}
H. Bateman,
Some recent researches on the motion of fluids, 
Monthly Weather Review, 43:4 (1915) 163-170. \\

\bibitem{49bochner}
S. Bochner,
Diffusion equations and stochastic processes, 
Proc Natl Acad Sci. U  S  A. Jul; 35, 7 (1949) 368-370.  \\

\bibitem{94boyadzhiev}
K. N. Boyadzhiev, 
Logarithms and imaginary powers of operators on Hilbert spaces,
Collect. Math. {\bf 45} 3 (1994) 287-300. \\

\bibitem{48burgers}
J. M. Burgers,
A mathematical model illustrating the theory of turbulence,
Adv. Appl. Mech.. 1 (1948) 171-199. \\

\bibitem{98cazenave}
T. Cazenave and A. Haraux,
An introduction to semilinear evolution equations,
Oxford University Press, 1998. \\

\bibitem{51cole}
J. D. Cole, 
On a quasi-linear parabolic equation occurring in aerodynamics,
Quart. Appl. Math. {\bf 9}  3 (1951) 225-236. \\

\bibitem{43dunford}
N. Dunford,
Spectral theory I, Convergence to projections,
Trans. Amer. Math. Soc. {\bf 54} (1943) 185-217.  \\

\bibitem{06forsyth}
A. R. Forsyth, 
Theory of differential equations. Part 4. 
Partial differential equations (Vol. 5-6), 1906. \\

\bibitem{03hasse}
M. Hasse,
Spectral properties of operator logarithms,
Math. Z. {\bf 245} 4 (2003) 761-779. \\

\bibitem{06hasse}
M. Hasse,
The functional calculus for sectorial operators,
Birkh\"auser, 2006. \\

\bibitem{74hirsch}
M. W. Hirsch and S. Smale,
Differential equations, dynamical systems, and linear algebra,
Academic Press, 1974. \\

\bibitem{50hopf}
E. Hopf,
The partial differential equation $u_t + u u_x = \mu u_{xx}$,
Comm. Pure and Appl. Math. {\bf 3} (1950) 201-230.  \\

\bibitem{17iwata-1}
Y. Iwata,
Infinitesimal generators of invertible evolution families,
Methods Func. Anal. Topology {\bf 23} 1 (2017) 26-36. \\

\bibitem{17iwata-3}
Y. Iwata,
Alternative infinitesimal generator of invertible evolution families, 
J. Appl. Math. Phys. {\bf 5} (2017) 822-830. \\

\bibitem{17iwata-2}
Y. Iwata,
Operator algebra as an application of logarithmic representation of infinitesimal generators
J. Phys.: Conf. Ser. 965 (2018) 012022. \\

\bibitem{18iwata-1}
Y. Iwata,
Abstract formulation of the Cole-Hopf transform,
Methods Funct. Anal. Topology {\bf 25} 2 (2019) 142-151. \\

\bibitem{18iwata-2}
Y. Iwata,
Relativistic formulation of abstract evolution equations,
AIP Conference Proceedings 2075  (2019) 100007. \\

\bibitem{19iwata-book}
Y. Iwata,
Operator topology,
A chapter of a book "Topology", IntechOpen, 2020 (DOI:10.5772/intechopen.92226).\\

\bibitem{19iwata}
Y. Iwata,
Unbounded formulation of the rotation group,
J. Phys. Conf. Ser. 1194 (2019) 012053.\\
\bibitem{20iwata}
Y. Iwata,
Abstract formulation of the Miura transform,
Mathematics, 2020, 8 747. \\

\bibitem{66kato}
T. Kato,
Perturbation Theory for Linear Operators,
Springer-Verlag, 1966. \\

\bibitem{70kato}
T. Kato,
Linear evolution equation of ''hyperbolic" type,
J. Fac. Sci. Univ. Tokyo {\bf 17} (1970) 241-258.  \\

\bibitem{73kato}
T. Kato,
Linear evolution equation of ''hyperbolic" type II,
J. Math. Soc. Japan {\bf 25} 4 (1973) 648-666.  \\

\bibitem{75kato}
T. Kato,
Quasi-linear equations of evolution, with applications to partial differential equations,
Lecture Notes in Math. 448, Springer-Verlag, 1975.\\

\bibitem{82kato}
T. Kato,
A short introduction to perturbation theory for linear operators,
Springer-Verlag, 1982. \\

\bibitem{72krein}
S. G. Krein,
Linear differential equations in Banach space (translated from Russian),. 
Transl. Math. Monogr. 29, Amer. Math. Soc., 1971. \\

\bibitem{11kreyszig}
E. Kreystig, H. Kreystig, and E. D. Norminton, Advanced engineering mathematics, John Wiley and Sons, 2011.
  
\bibitem{82landau}
L. D. Landau and E. M. Lifshitz,
Mechanics Third Edition: Volume 1 (Course of Theoretical Physics) 
Elsevier, 1982. \\

\bibitem{73sagle}
A. A. Sagle and R. E. Walde, 
Introduction to Lie Groups and Lie algebras,
Academic Press, 1973. \\

\bibitem{61lions}
J. L. Lions, 
Equations differentielles operationnelles et problemes aux limits,
Springer-Verlag, 1961.  \\

\bibitem{72lions}
J. L. Lions and E. Magenes, 
Non-homogeneous boundary value problems and applications,
Springer-Verlag, 1972.  \\


\bibitem{61lumer}
G. Lumer and R. S. Phillips,
Dissipative operators in a Banach spaces,
Pacific J. Math. {\bf 11} (1961) 679-698. \\

\bibitem{01martinez}
C. Martinez and M. Sanz,
The theory of fractional powers of operators,
North-Holland, 2001. \\

\bibitem{85mizohata}
S. Mizohata,
On the Cuachy problem,
Academic Press, 1985. \\

\bibitem{69nollau}
V. Nollau,
\"Uber den Logarithmus abgeschlossener Operatoren in Banachschen R\"aumen, 
Acta Sci. Math. {\bf 300} (1969) 161-174.  \\

\bibitem{00okazawa-1}
N. Okazawa,
Logarithms and imaginary powers of closed linear operators,
Integral Equations and Operator Theory {\bf 38} 4 (2000) 458-500. \\

\bibitem{00okazawa-2}
N. Okazawa,
Logarithmic characterization of bounded imaginary powers,
Progress in Nonlinear Differential Equations and Their Applications, {\bf 42} (2000) 229-237.  \\

\bibitem{83pazy}
A. Pazy, 
Semigroups of linear operators and application to partial differential equations, 
Springer-Verlag, 1983. \\

\bibitem{52phillips}
R. S. Phillips,
On the generation of semi-groups of linear operators,
Pacific J. Math. 2 (1952) 343-369. \\

\bibitem{01pruess}
J. Pr\"uss and R. J. Schnaubelt,
Solvability and Maximal Regularity of Parabolic Evolution Equations with Coefficients Continuous in Time,
Math. Anal. Appl. {\bf 256} (2001) 405-430. \\

\bibitem{79tanabe}
H. Tanabe,
Equations of evolution,
Pitman, 1979. \\

 \bibitem{51taylor}
A. E. Taylor, 
Spectral theory of closed distributive operators.
Acta Math. {\bf 84} (1951) 189-224.  \\

\bibitem{97temam}
R. Temam,
Infinite-dimensional dynamical systems in mechanics and physics, 2nd edition,
Springer-Verlag, 1997. \\

\bibitem{50weyl}
H. Weyl,
The theory of groups and quantum mechanics,
Dover publications, 1950. \\

\bibitem{57yamanouchi}
T. Yamanouchi,
The rotation group and its representation (in Japanese),
Iwanami-Shoten publishers: 1957. \\

\bibitem{60yosida}
K. Yosida
Fractional Powers of Infinitesimal Generators and the Analyticity of the Semi-groups Generated by Them, 
Proceedings of the Japan Academy, {\bf 36} 3 (1960) 86-89. \\


\bibitem{65yosida}
K. Yosida,
Functional Analysis, 
Springer-Verlag, 1965. \\

\end{thebibliography}
\end{document}